\RequirePackage{fix-cm}
\documentclass[smallextended]{svjour3}       

\smartqed  

\usepackage{graphicx}

\usepackage[fleqn]{amsmath}
\usepackage{amsfonts}

\usepackage[section]{placeins} 

\usepackage{algorithmic}
\usepackage{algorithm}

\newcommand{\algorithmicbreak}{\textbf{break}}
\newcommand{\BREAK}{\STATE \algorithmicbreak}
\algsetup{linenodelimiter=.}

\usepackage{url}
\usepackage[colorlinks=true,linkcolor=blue,citecolor=blue,pdfpagelabels=true]{hyperref}
\usepackage{doi}

\newcommand{\R}{\mathbb{R}}
\newcommand{\diff}[1]{\mathrm{d}#1}
\newcommand{\dt}{\mathrm{d}t}

\newcommand{\norm}[1]{\lVert #1 \rVert}

\renewcommand{\exp}[1]{\mathrm{e}^{#1}} 
\newcommand{\T}{^{\mathsf{T}}}  

\newcommand{\cF}{\mathcal{F}}
\newcommand{\cG}{\mathcal{G}}

\newcommand{\cS}{\mathcal{S}}
\newcommand{\cStr}{\mathcal{S}_{\text{Strang}}}
\newcommand{\cL}{\mathcal{S}_{\text{Lie}}}
\newcommand{\cSym}{\mathcal{S}_{\text{sym}}}
\newcommand{\cAsym}{\mathcal{S}_{\text{asym}}}

\newcommand{\cTF}{\mathcal{T}_{\cF}}
\newcommand{\cTG}{\mathcal{T}_{\cG}}
\newcommand{\cTFG}{\mathcal{T}_{\cF + \cG}}

\spdefaulttheorem{assumption}{Assumption}{\bf}{\rm}

\DeclareMathOperator{\blkdiag}{blkdiag}
\DeclareMathOperator{\trace}{trace}

\begin{document}

\title{Adaptive high-order splitting schemes for large-scale differential Riccati equations}

\author{Tony Stillfjord}

\institute{T. Stillfjord \at
              Mathematical Sciences, Chalmers University of Technology and the University of Gothenburg, SE-412 96 G\"{o}teborg, Sweden \\
              Tel.: +46-31-7725304\\
              Fax: +46-31-161973\\
              \email{tony.stillfjord@gu.se}
}

\date{Received: date / Accepted: date}

\maketitle

\begin{abstract}
We consider high-order splitting schemes for large-scale differential Riccati equations. Such equations arise in many different areas and are especially important within the field of optimal control. In the large-scale case, it is critical to employ structural properties of the matrix-valued solution, or the computational cost and storage requirements become infeasible. Our main contribution is therefore to formulate these high-order splitting schemes in a efficient way by utilizing a low-rank factorization. Previous results indicated that this was impossible for methods of order higher than 2, but our new approach overcomes these difficulties. In addition, we demonstrate that the proposed methods contain natural embedded error estimates. These may be used e.g.\ for time step adaptivity, and our numerical experiments in this direction show promising results.
\end{abstract}





\section{Introduction} \label{sec:introduction}
We consider differential Riccati equations (DREs) of the form
\begin{equation}
  \label{eq:DRE}
  \dot{P} = A^T P + PA + Q - P S P, \quad P(0) = P_0,
\end{equation}
where the solution $P(t)$ is matrix-valued and $A$, $Q$ and $S$ are given matrices. Such equations arise in many different areas, e.g.\ in 
optimal/robust control, optimal filtering, spectral factorizations, $\mathbf{H}_\infty$-control, differential games, etc.~\cite{AbouKandil_etal2003,BasarBernhard1995,IchikawaKatayama1999,Petersen_etal2000}.

 A typical application is a linear quadratic regulator (LQR) problem, where one seeks to control the output $y = Cx$ given the state equation $\dot{x} = Ax + Bu$ by varying the input $u$. In the case of a finite time cost function, 
\begin{equation*}
  J(u) = \int_{0}^{T}{x(t)^T R_x x(t) + u(t)^T R_u u(t) \dt} , 
\end{equation*}
where $R_x$ and $R_u$ are given matrices, it is well known that the optimal input $u^*$ is given in state feedback form. In particular, $u^*(t) = -R_u^{-1}B^TP(T-t)x(t)$, where $P$ is the solution to the DRE~\eqref{eq:DRE} with the specific matrices $Q = C^TR_xC$ and $S = B R_u^{-1}B^T$. We note that the situation $M\dot{x} = Ax + Bu$ can be handled in a straightforward way without explicitly inverting the mass matrix $M$, see e.g.~\cite{Stillfjord2015}.

In this paper, we are interested in the large-scale setting. Even if ${A \in \R^{N \times N}}$ is sparse, the solution $P$ is typically dense. Hence, a ``large'' dimension $N$ is here considerably smaller than the number of components which would be considered large for a vector-valued ODE. A naive method that works well for the small-scale case would run into storage problems already for  $N = 10000$ and be computationally expensive long before that. Recently, many non-naive methods have been proposed for DREs and similar problems, e.g.\ matrix-valued BDF and Rosenbrock methods~\cite{BennerMena2016,BennerMena2013}, splitting schemes~\cite{Stillfjord2015,Mena_etal2017} and Krylov projection methods~\cite{KoskelaMena2017,Guldogan_etal2016}. 
The latter are a generalization of the Krylov approach to algebraic Riccati equations and Lyapunov equations~\cite{DruskinKnizhermanSimoncini2011,HeyouniJbilou2009,SimonciniSzyldMonsalve2014}. Other methods for such equations, like invariant subspace techniques~\cite{AmodeiBuchot2010,BennerBujanovic2016,YidingSimoncini2015}, typically also generalize to the DRE case by using time-stepping methods of either one- or multi-step type.
Further useful references may be found in the recent surveys~\cite{BennerSaak2013,Simoncini2016}. In general, all these methods rely on the fact that the dense solution possesses certain structure. In particular, the solution is positive semi-definite, and in all practical applications it also has low rank. This allows us to factorize $P = ZZ^T$ where $Z$ is a matrix with many fewer columns than $P$. A main idea in all the algorithms listed above is then to only do computations on the factor $Z$ and never actually form the product $ZZ^T$.

Further, we are interested in different types of splitting schemes, since the equation has a natural division into two parts:
\begin{equation*}
\dot{P} = \cF P + \cG P, \quad \text{where} \quad \cF P = A^T P + PA + Q \quad \text{and} \quad \cG P = - PSP.
\end{equation*}
While the full problem is rather difficult, the subproblems
\begin{align}
\label{eq:subF}  \dot{P} &= \cF P,  \quad P(0) = P_0, \quad \text{and} \\
\label{eq:subG}  \dot{P} &= \cG P,  \quad P(0) = P_0,
\end{align}
are separately much easier and cheaper to solve. In fact, as demonstrated in~\cite{Stillfjord2015} there exist closed-form expressions for the solutions to both subproblems that are amenable to low-rank computations. In the following, we will denote the solution operator to the full problem by $\cTFG$ and to the subproblems by $\cTF$ and $\cTG$; thus for example the solution to~\eqref{eq:subF} at time $t$ is given by $\cTF(t) P_0$.

To introduce the simplest splitting schemes and our notation, we first discretize the time interval $[0, T]$ by $n$ equidistant time steps of size $h$ and set $t_j = j h$. Then the approximation to $\cTFG(t_j)P_0$ by the Lie splitting scheme is given by $\cL(h)^j P_0$, where
\begin{equation*}
  \cL(h) = \cTF(h) \cTG(h).
\end{equation*}
That is, we switch back and forth between the affine subproblem and the nonlinear subproblem. A more accurate approximation is given by the Strang splitting scheme, defined by the time stepping operator
\begin{equation*}
  \cStr(h) = \cTG(h/2) \cTF(h) \cTG(h/2).
\end{equation*}
In both cases, we may interchange the order of the $\cF$ and $\cG$ operators.
For a more thorough introduction to splitting schemes in general, we refer to~\cite{HundsdorferVerwer2003}.

It can be shown as in~\cite{HundsdorferVerwer2003} that the Lie splitting is first-order convergent and the Strang splitting second-order convergent, i.e.\ the errors satisfy
\begin{equation*}
  \norm{\cL(h)^j P_0 - \cTFG(t_j)P_0} \le C h \quad \text{and} \quad \norm{\cStr(h)^j P_0 - \cTFG(t_j)P_0} \le C h^2.
\end{equation*}
In general one can also consider higher-order schemes, but so far this has not been done for DREs. This is due to the fact that multiplicative splitting schemes of the form $\cTF(\alpha_1 h) \cTG(\beta_1 h) \cdots \cTF(\alpha_s h) \cTG(\beta_s h) $ require that some coefficients $\alpha_j$, $\beta_j$ are either negative or complex~\cite{BlanesCasas2005,HansenOstermann2009_2}, which is not compatible with the low-rank implementation.

The first main contribution of this work is therefore to demonstrate that a new type of additive splitting schemes introduced in~\cite{DeLeo_etal2016} allows for arbitrary high order schemes to be implemented efficiently in a low-rank DRE setting. These schemes are of the form
\begin{equation*}
  \gamma_1 \cTF(h) \cTG(h) + \gamma_2 \big(\cTF(h/2) \cTG(h/2)\big)^2 + \cdots \gamma_s \big(\cTF(h/s) \cTG(h/s)\big)^s
\end{equation*}
and thus only utilize positive step sizes. A minor drawback is that the approximations are no longer guaranteed to be positive semi-definite, since the coefficients $\gamma_j$ may be negative. This prohibits the use of a $ZZ^T$-factorization, and we therefore outline the changes necessary to instead consider a so-called $LDL^T$-factorization (cf.~\cite{LangMenaSaak2015}). 

The second main contribution lies in the observation that these splitting schemes contain natural lower-order embedded methods, which allows for cheap and easy error estimation. We utilize this to construct high-order splitting schemes with adaptive time-stepping, i.e.\ the time steps $h_j = t_{j+1} - t_j$ are no longer equidistant but chosen as large as possible while keeping the error below a given tolerance. Modifying the step size can greatly increase the efficiency, but only if the computational cost of changing the step size is small. We therefore outline which quantities can be precomputed or recomputed cheaply, and describe efficient updating strategies for the quantities that necessarily change with each step.

A brief outline of the paper is as follows. In Section~\ref{sec:lowrank} we state the basic assumptions on the given data and review the use of the $ZZ^T$- and $LDL^T$-factorizations for low-order splitting schemes. The issues that arise when considering higher-order multiplicative splitting schemes are outlined in Section~\ref{sec:splitting}, wherein we also present the new type of additive schemes that eliminate these issues. Error estimates and different kinds of time step adaptivity are discussed in Section~\ref{sec:adaptivity} and an algorithm summarizing the complete implementation is presented. In Section~\ref{sec:experiments}, several numerical experiments demonstrate the validity of the implementation, the efficiency of the methods and the use of adaptive time-stepping. Finally, we collect some conclusions in Section~\ref{sec:conclusions}.

\section{Low-rank factorizations} \label{sec:lowrank}
The first assumption we make on the problem data is the following:
\begin{assumption}\label{ass:main}
  The matrices $A$, $Q$, $S$ and the initial condition $P_0$ all belong to $\R^{N \times N}$. In addition, $Q$, $S$ and $P_0$ are symmetric and positive semi-definite.
\end{assumption}
This implies the existence and uniqueness of a solution $P$ to the DRE~\eqref{eq:DRE} such that $P(t)$ is also symmetric and positive semi-definite for all $t \ge 0$~\cite[Theorem 4.1.6]{AbouKandil_etal2003}. An important example of when Assumption~\ref{ass:main} is satisfied is the LQR setting from the introduction, with $R_x$ and $R_u$ both symmetric positive definite.
Secondly, we assume that the solution has the low-rank property:
\begin{assumption}\label{ass:lowrank}
  For each $t \in [0, T]$, the rank of the solution $P(t)$ is at most $r \ll N$ and the rank of $Q$ is $r_Q \ll N$.
\end{assumption}
To the author's knowledge there are currently no known useful criteria on the data in the DRE setting which guarantee that Assumption~\ref{ass:lowrank} is fulfilled. However, such low-rank structure is observed in all practical applications, e.g.\ in LQR problems where $B \in \R^{N \times m_B}$ and $C \in \R^{m_C \times N}$ with $m_B, m_C \ll N$. Recently, some results in this direction has been established for algebraic Riccati equations, i.e.\ the stationary version of~\eqref{eq:DRE}, in~\cite{BennerBujanovic2016}. These are generalizations of results for Lyapunov equations~\cite{AntoulasSorensenZhou2002,SorensenZhou2002} and it seems likely that further generalizations to the DRE setting could be made.

Assumptions~\ref{ass:main} and~\ref{ass:lowrank} imply that we can low-rank factorize $P(t) = Z(t)Z(t)^T$ and $Q = qq^T$ with $Z(t) \in \R^{N \times r}$ and $q \in \R^{N \times r_Q}$. Similarly, as demonstrated in~\cite{Stillfjord2015} we can low-rank factorize also the approximations $\cL(h) P_0$ and $\cStr(h) P_0$. This is based on factorizing the exact solutions to the subproblems~\eqref{eq:subF}-\eqref{eq:subG}, for which we have the closed-form expressions
\begin{align*}
  \cTF(h)P_0 &=  \exp{h A\T }P_0 \exp{hA}  + \int_{0}^{h}{\exp{s A\T } Q \exp{sA} \diff{s}} \quad \text{and}\\
  \cTG(h)P_0 &=  (I + hP_0S)^{-1}P_0.
\end{align*}
The latter expression quickly yields an explicit factorization while the former requires that the integral is approximated by a quadrature formula, whereafter column compression is applied.

Considering instead a so-called $LDL^T$-factorization where $L(t) \in \R^{N \times r}$ and $D(t) \in \R^{r \times r}$ is beneficial for many schemes~\cite{LangMenaSaak2015}, because it can decrease the amount of computations. This is true also for splitting schemes. Assuming that $P_0 = LDL^T$ and considering first the nonlinear subproblem, we have
\begin{equation*}
  \cTG(h)P_0 = (I + hLDL^TS)^{-1}LDL^T    = L (I + hDL^TSL)^{-1}DL^T,
\end{equation*}
by use of a simplified version of the Woodbury matrix inversion formula~\cite{Hager1989}. Thus $\hat{L}\hat{D}\hat{L}^T$ is a low-rank factorization of the solution to the nonlinear subproblem, where $\hat{L} = L$ and $\hat{D} = (I + hDL^TSL)^{-1}D$. In contrast to the $ZZ^T$ situation, there exist matrices $D$ and $L$ such that $I + hDL^TSL$ is not invertible for all $h$. However, it certainly is for all $h < 1/\rho(DL^TSL)$, where $\rho$ denotes the spectral radius, and therefore the step size can always be chosen such that $\hat{D}$ is well defined. We note that even for large time steps, this theoretical issue has not yet been observed in practice. We also note that this formulation is cheaper to compute than the corresponding $ZZ^T$-factorization, since it is no longer necessary to compute a Cholesky factorization of the inverse.

Considering next the affine subproblem and assuming that $Q = L_Q D_Q L_Q^T$, we have
\begin{align*}
  \cTF(h)P_0 &=  \exp{h A\T } LDL^T \exp{hA}  + \int_{0}^{h}{\exp{s A\T } L_Q D_Q L_Q^T \exp{sA} \diff{s}} \\
             &= L_1 D L_1^T + \int_{0}^{h}{L(s) D_Q L(s)^T \diff{s}} \\
             &\approx  L_1 D L_1^T + \sum_{k=1}^{n_Q}{ w_k L(s_k) D_Q L(s_k)^T} ,
\end{align*}
where $L_1 = \exp{h A\T } L$, $L(s) = \exp{s A\T } L_Q$ and $(s_k, w_k)$ are the $n_Q$ nodes and weights of a quadrature formula. We choose the parameters such that the error in this approximation is negligible with respect to the splitting error; for a splitting scheme of order $p$ we typically choose a quadrature formula of order $p+1$. 
For efficiency, the structure of $A$ (sparsity, bandedness, etc.) should be taken into account when computing the terms $L_1$ and $L(s)$. In our tests, we simply use a 5th-order implicit Runge-Kutta method with a crude error estimate based on halving the internal step size. It seems likely, however, that an approach based on e.g. Krylov subspaces or the Leja point method (see e.g.~\cite{Caliari_etal2014}) would be even more efficient, especially if subspaces from previous steps can be (partially) reused. We note that these terms do not need to be computed to full precision, but like for the integral term their errors should be negligible in comparison to the splitting error.
Then, similarly to the $ZZ^T$-case, setting
\begin{equation*}
  \tilde{L} =
  \begin{bmatrix}
    L_1 & L(s_1) & \cdots & L(s_{n_Q})
  \end{bmatrix}
\quad \text{and} \quad
\tilde{D} =
\begin{bmatrix}
  D &        &         & \\
    & w_1D_Q  &         & \\
    &        &  \ddots & \\
    &        &         & w_{n_Q} D_Q
\end{bmatrix}
\end{equation*}
means that $\tilde{L}\tilde{D}\tilde{L}^T$ is  a low-rank approximation of the solution to the affine subproblem. After forming $\tilde{L}$ and $\tilde{D}$, column-compression should be applied to eliminate any unnecessary columns. We refer to~\cite{LangMenaSaak2015} for an efficient way to do this.

\section{High-order splitting schemes} \label{sec:splitting}
Let us now consider low-rank factorization of higher-order multiplicative splitting schemes like the Lie and Strang splitting schemes. Let 
  \begin{equation*}
    \cS(h) = \cTF(\alpha_1 h) \cTG(\beta_1 h) \cdots \cTF(\alpha_s h) \cTG(\beta_s h)
  \end{equation*}
with $s$ and the coefficients $\{\alpha_k\}_{k=1}^{s}$, $\{\beta_k\}_{k=1}^{s}$ chosen such that $\cS(h)$ is a splitting scheme of order $p \ge 3$. Then the coefficients must include either negative or complex values~\cite{BlanesCasas2005,HansenOstermann2009_2}. In the first case, computing $\exp{\gamma h A^T}P_0$ for such a negative coefficient $\gamma$ corresponds to taking a negative time-step for the system $\dot{x} = A^Tx$. If $A$ e.g.\ corresponds to a discretization of the Laplacian (a common application) we are thus solving the heat equation backwards in time, which is ill-posed. It is therefore only possible to consider the class of problems where $A$ corresponds to the discretization of an analytic operator, but even in this case the evaluation of $(I + hZ^TSZ)^{-1}$ or $(I + hDL^TSL)^{-1}$ tends to yield step size restrictions. We therefore do not think that this is a worthwhile direction of research to pursue.

In the case of a $ZZ^T$-factorization, a complex coefficient $\gamma$ destroys the structure of $I + \gamma hZ^TSZ$ and we can only factorize it in very special cases. Considering instead an $LDL^T$-factorization leads to problems with complex arithmetic: If $L$ and $D$ are real, the approximation $\hat{L} \hat{D}\hat{L}^T$ to $\cTG(\gamma h)LDL^T$ will have $\hat{L}$ real but $\hat{D}$ complex-valued. Such input to the affine subproblem will then lead to both $\hat{L}$ and $\hat{D}$ being complex-valued. Once this is the case, we not only have to do computations fully in complex arithmetic but we also have issues with column compression since the complex values do not match the ``transpose''-formulation. Switching instead to a complex $LDL^H$-factorization results in similar issues. Like negative coefficients, using complex coefficients thus does not seem worthwhile.

However, the necessity of negative or complex coefficients only holds for the type of multiplicative splitting schemes mentioned above. Recently, a new type of additive splitting schemes was introduced in~\cite{DeLeo_etal2016}.
These are either of the asymmetric type
\begin{equation} \label{eq:asym}
  \cAsym^s(h) = \sum_{k = 1}^{s}{ \gamma_k \big(\cTF(h/k) \cTG(h/k) \big)^k},
\end{equation}
which are of order $s$ if the coefficients $\gamma_1, \ldots, \gamma_s$ are chosen appropriately, and the symmetric type
\begin{equation} \label{eq:sym}
  \cSym^{2s}(h) = \sum_{k = 1}^{s}{ \gamma_k \bigg( \big(\cTF(h/k) \cTG(h/k) \big)^k + \big(\cTG(h/k) \cTF(h/k) \big)^k \bigg)},
\end{equation}
which are of order $2s$. (We only consider the case of minimal number of stages here. One might of course add extra stages in order to improve the local error structure, but given the form of the schemes it would then make more sense to instead increase the order.) In both cases, the roles of $\cF$ and $\cG$ may be interchanged.

 At first sight these methods may look computationally expensive. However, (as noted in~\cite{DeLeo_etal2016}) if we have the possibility to work in parallel then taking one step with either method is only as expensive as taking $s$ Lie splitting steps. More important is that they only require real, positive step sizes. This eliminates all the issues listed above, and allows us to consider splitting schemes for DREs of arbitrarily high order.

Because the coefficients $\{\gamma_k\}_{k=1}^s$ may include negative values, using a $ZZ^T$-factorization to formulate these methods is impossible. However, instead using an $LDL^T$-factorization is not only possible but rather straightforward after the preliminary work in the previous section. The only additional computational work is a column compression step after forming the linear combinations. In an optimized code, most of this work could additionally be done while waiting for the slowest processor that takes $s$ steps to finish. Using a higher-order method also requires us to compute terms of the form $\exp{\gamma h A^T}L$ more accurately (unless we also increase the step size $h$ and thereby the error), and to use a higher-order quadrature formula to approximate the integral term in $\cTF(h)P_0$.

We also note here that while the $LDL^T$-factorization does not guarantee that the approximations are positive semi-definite, in practice this still seems to hold. This is likely due to the fact that the approximations are very close to the solution of the full problem, which is guaranteed to be positive semi-definite.

\section{Time adaptivity} \label{sec:adaptivity}
An additional major feature of the schemes~\eqref{eq:asym} and~\eqref{eq:sym} is the existence of natural embedded lower-order methods. This seems to have been overlooked by~\cite{DeLeo_etal2016}. For example, the scheme
\begin{equation*}
  \cAsym^2(h) = -1 \big(\cTF(h) \cTG(h) \big) + 2 \big(\cTF(h/2) \cTG(h/2) \big)^2
\end{equation*}
is of order $2$, and it obviously contains the first-order method $\cTF(h) \cTG(h)$. This holds true for all the schemes, symmetric and asymmetric. In general, neglecting the last terms of the sum and using other coefficients $\{\gamma_k\}$ yields embedded methods of order $p-1$ in the asymmetric case, and of order $2p-2$ in the symmetric case with $p = 2, 3, \ldots, s$. Since these lower-order approximations are simply linear combinations of previously computed terms they are cheap to compute. In our case, the only extra computational effort is a column compression step.

The embedded methods yield natural error estimates. For example, we have
\begin{align*}
 & \big(\cAsym^s(h)P_0 - \cTFG(h)P_0 \big) - \big(\cAsym^{s-1}(h)P_0 - \cTFG(h)P_0 \big) \\
  &\qquad  = \Phi^s(P_0) h^{s+1} + \mathcal{O}(h^{s+2}) - \Phi^{s-1}(P_0) h^{s} + \mathcal{O}(h^{s+1}) \\
&\qquad = - \Phi^{s-1}(P_0) h^{s} + \mathcal{O}(h^{s+1}),
\end{align*}
where $\Phi^s$ and $\Phi^{s-1}$ are the principal error functions of the two methods. Thus the difference $\cAsym^s(h)P_0 - \cAsym^{s-1}(h)P_0$ is a local error estimate of order $s-1$. In the symmetric case we instead get an error estimate of order $2s-2$.

These error estimators may be used to control the size of $h$, with the aim of keeping the local error below a certain tolerance while minimizing the computational effort. There are many different kinds of such controllers, see e.g.~\cite{Gustafsson_etal1988,Soderlind2002}. As an example, we choose a simple PI-controller which typically provides a smoother step size sequence than the commonly used deadbeat I-controller. It is given by~\cite{Gustafsson_etal1988,Soderlind2002}
\begin{equation*}
  h_{n+1} = \bigg( \frac{\epsilon \, \text{TOL}}{e_{n+1}}\bigg)^{k_I} \bigg( \frac{e_n}{e_{n+1}}\bigg)^{k_P} h_n,
\end{equation*}
where $h_n$ is the $n$:th time step, $e_n$ is the error estimate at $t_n$, $\text{TOL}$ is the desired accuracy (tolerance) and $\epsilon$ is a safety factor . The parameters $(k_I, k_P)$ determine the characteristics of the controller such as responsiveness and robustness. In our numerical experiments we set $\epsilon = 0.9$ and $(k_I, k_P) = (0.2/p, 0.2/p)$, where $p$ is the order of the error estimate. These are similar to the values recommended for explicit Runge-Kutta methods when using the error per unit step strategy~\cite{Gustafsson1991}. Clearly, these are not optimal values for splitting schemes, but an in-depth investigation for a variety of typical problems is out of the scope of this paper.

The evaluation of $\cTG(h)P_0$ requires the same effort whether the step size is varying or not. Evaluating $\cTF(h)P_0$, on the other hand, requires that the approximation of the integral term
\begin{equation*}
  I_Q(h_n) = \int_{0}^{h_n}{\exp{s A\T } Q \exp{sA} \diff{s}}
\end{equation*}
is recomputed in every step, while it previously could be precomputed. We note that typically the rank of $Q$ is sufficiently small in relation to the rank of the solution approximation that this extra computational cost is small and easily outweighed by the benefits of adaptivity. Nevertheless, we suggest here a strategy to decrease this cost further.

We need to compute $\sum_{k=1}^{n_Q}{ w_k L(s_k) D_Q L(s_k)^T}$ where $L(s) = \exp{s A\T } L_Q$ for given nodes $s_k$ and weights $w_k$. The main idea now is to change only a few nodes (and thereby also the weights) in each step, such that the interval $[0, h_n]$ is covered as evenly as possibly. For a quadrature rule of order $p$ we need $n_Q = p+1$ nodes if we do not place the nodes optimally, in contrast to e.g. Gaussian quadrature which would need roughly half as many. However, by storing the computed matrices $L(s_k)$ and keeping most of the nodes unchanged we will still decrease the overall computation cost. Thus, we define the initial nodes by $s_k = \frac{k h_1}{p}$ for $k = 0, \ldots, p$ and compute the initial weights from
\begin{equation} \label{eq:quad_weights}
  \begin{bmatrix}
    s_0^0   & s_1^0 & \cdots & s_p^0 \\
    s_0^1   & s_1^1 &        & \vdots     \\
    \vdots &       & \ddots        &      \\
    s_0^p   & \cdots      &         & s_p^p
  \end{bmatrix}
  \begin{bmatrix}
    w_0\\
    w_1\\
    \vdots\\
    w_p
  \end{bmatrix}
  =
  \begin{bmatrix}
    h_n\\
    h_n^2 / 2\\
    \vdots\\
    h_n^{p+1} / (p+1)
  \end{bmatrix}
  \end{equation}
with $n = 1$.
Then, to update these nodes and weights given a new $h_n$, we follow the procedure outlined in algorithmic form in Algorithm~\ref{alg:update_IQh}. (In Algorithm~\ref{alg:update_IQh} and in the following, $\blkdiag$ denotes the block diagonal operator, i.e.\ it places its block arguments on the diagonal of an otherwise zero matrix.)

Essentially, we add a node at $h_n$ if the interval increases, and then remove the node which makes the remaining sequence as close to equidistributed as possible. Similarly, if the interval decreases, we iteratively relocate the nodes that are outside the new interval to the midpoints of the largest gaps between the nodes in the new interval. In order to ensure that the nodes $s_k$ cover the interval $[0, h_n]$ well, we recompute the whole sequence if the step size changes by more than $25 \%$. 
We note that we could of course, in theory, store all the previously computed $L(s_k)$ and use increasingly high-order quadrature formulae. However, this would yield a major increase in the storage requirements while having little effect on the overall accuracy.

\begin{algorithm}
  \caption{Updating the low-rank factorization of $I_Q(h_n)$}
  \label{alg:update_IQh}
  \begin{algorithmic}[1]
    \REQUIRE Old and new step sizes $(h_{n-1}, h_n)$, previous nodes $\{s_k\}_{k=0}^{p}$, matrices $\{L(s_k)\}_{k=0}^{p}$
    \IF{$h_n \le 0.8 \, h_{n-1}$ \OR $h_n \ge 1.25\, h_{n-1}$ }
    \STATE Set $\hat{s}_k = \frac{k h_n}{p}$, $k=0, \ldots, p$
    \STATE Recompute all $L(\hat{s}_k)$
    \ELSIF{$h_n > h_{n-1}$}
    \STATE Set $\hat{s}_k = s_k$, $k = 0,\ldots, p$, and $\hat{s}_{p+1} = h_n$
    \STATE Remove the node $\hat{s}_j$ such that $\displaystyle d_j = \min_{0 \le k \le p+1}{d_k}$, where $d_0 = \hat{s}_1$, $d_{p+1} = h_n - \hat{s}_p$
 and $d_k = \hat{s}_{k+1} - \hat{s}_{k-1}$ for $k=1,\ldots,p$
    \IF{the new node $\hat{s}_{p+1}$ was removed}
    \STATE Set $L(\hat{s}_k) = L(s_k)$, $k = 0, \ldots, p$
    \ELSE
    \STATE Compute $L(h_n) = \exp{h_n A^T}L_Q$
    \STATE Set $L(\hat{s}_k)$ to the matrices $L(s_k)$ and $L(h_n)$ that match the nodes $\hat{s}_k$
    \ENDIF
    \ELSIF{$h_n < h_{n-1}$}
    \STATE Find the number $\tilde{n}$ of nodes to recompute: $\tilde{n} = p + 1 - j$, with $j = \max \{k : s_k \le h_n \}$
    \IF{$\tilde{n} = 0$}
    \STATE Set $\hat{s}_k = s_k$ and $L(\hat{s}_k) = L(s_k)$ for $k = 0, \ldots, p$, i.e.\ do nothing
    \ELSE
    \STATE Set $\hat{s}_k = s_k$ and $L(\hat{s}_k) = L(s_k)$ for $k = 0, \ldots, j$  
    \FOR{$l = 1, \ldots, \tilde{n}$ }
    \STATE Find $i$ such that $\displaystyle d_i = \max_{0 \le k \le j+l}{d_k}$, where $d_0 = \hat{s}_0$, $d_{j+l} = h_n - \hat{s}_{j+l-1}$ and $d_k = \hat{s}_{k} - \hat{s}_{k-1}$ for $k=1,\ldots,j+l-1$
    \STATE Add a new node $\hat{s}_{j+l}$ at $\hat{s}_0/2$ if $i = 0$, at $(h_n + \hat{s}_{j+l-1})/2$ if $i = j+l$ or  at $(s_{i} + s_{i-1})/2$ if $1 \le i \le j+l-1$
    \STATE Compute $L(\hat{s}_{j+l}) = \exp{\hat{s}_{j+l} A^T}L_Q$
    \STATE Reorder $\hat{s}_k$ and $L(\hat{s}_k)$ so that the nodes are increasing

    \ENDFOR
    \ENDIF
    
    \ENDIF
   
    \STATE Compute new weights $\{\hat{w}_k\}_{k=0}^{p}$ from Equation~\eqref{eq:quad_weights}
    \STATE Form $\hat{L} = \begin{bmatrix} L(\hat{s_0}) & \cdots &  L(\hat{s_p}) \end{bmatrix} $ and $\hat{D} = \blkdiag(\hat{w}_0 D_Q, \ldots, \hat{w}_p D_Q)$
    \STATE Column-compress $\hat{L}$ and $\hat{D}$
    \ENSURE New nodes $\{\hat{s}_k\}_{k=0}^{p}$, matrices $\{L(\hat{s}_k)\}_{k=0}^{p}$, weights $\{\hat{w}_k\}_{k=0}^{p}$, matrices $\hat{L}$ and $\hat{D}$ such that $\hat{L}\hat{D} \hat{L}^T \approx I_Q(h_n)$
  \end{algorithmic}
\end{algorithm}

\begin{remark}
  We note that in e.g.\ a real-world optimal control problem, it is frequently the case that the state of the system is sampled at regular, predetermined intervals. The feedback control thus needs the solution of the corresponding DRE at these specific times. This suggests that a constant, matching step size should be employed, or that the adaptive step size is restricted. Neither approach is desirable; the former is inefficient compared to the adaptive approach, and the latter destroys the smooth time step sequence the PI controller is intended to provide. However, assuming that the exact solution is sufficiently regular, we may still use the more efficient adaptive time stepping and simply interpolate the computed approximations to find the values at the desired times. For example, assume that we use piecewise linear interpolation. Then on the interval $[t_{n-1}, t_n]$, the error between the interpolant $P_I$ and the exact solution $P$ is bounded by
\begin{equation*}
  \norm{P_I(t) - P(t)} \le \text{TOL } + h_n^2 \sup_{s \in [t_{n-1}, t_n]}{\norm{\ddot{P}(s)}}/8.
\end{equation*}
To estimate the second term,  we can first use the available approximation $P_n \approx P(t_n)$ to estimate $\dot{P}(t_n) \approx A^T P_n + P_n A + Q - P_n S P_n$. Then by differentiating Equation~\eqref{eq:DRE} we get $\ddot{P} = A^T\dot{P} + \dot{P}A - \dot{P} S P - P S \dot{P}$, from which we can estimate $\ddot{P}(t_n)$. Both of these operations may be low-rank factorized; if $P_n = LDL^T$ then $\dot{P}(t_n) \approx \tilde{L}\tilde{D}\tilde{L}^T$ where 
$
\tilde{L} =
\begin{bmatrix}
  A^TL\; & L\; & L_Q
\end{bmatrix}
$, and $\ddot{P}(t_n) \approx \hat{L}\hat{D}\hat{L}^T$ where 
$
\hat{L} =
\begin{bmatrix}
  A^T\tilde{L}\; & \tilde{L}\; & L
\end{bmatrix}
$
.
 Thus the norm of $\ddot{P}$ may be estimated efficiently by two applications of $A^T$ and two column compression operations. This estimation may be incorporated into the step-size controller to automatically ensure that the interpolation error is bounded by a fixed tolerance. Whether this is cost-effective or not is of course heavily dependent on the rank of the approximation, and thus of the problem data.

To actually compute the interpolant in a low-rank setting, we note that if $P_{n-1} = L_{n-1} D_{n-1} L_{n-1}^T$ and $P_n = L_n D_n L_n^T$ then $L = \begin{bmatrix}  L_{n-1} & L_n \end{bmatrix}$ and $D = \begin{bmatrix} \alpha D_{n-1} & 0 \\ 0 & (1-\alpha)D_{n} \end{bmatrix}$ constitute a low-rank factorization of $\alpha P_{n-1} + (1-\alpha) P_n$, so that this computation comes at the cost of one column compression step. This interpolation procedure is less straightforward if the setting is generalized to that of time-varying matrices. However, in that case the strategy of sampling the system at constant time intervals is also rather dubious.
\end{remark}
\begin{remark}
  It is enough if the terms involved in one step of the splitting methods are of the same accuracy as the local error. Therefore, the error estimates may additionally be used to determine the optimal tolerances for column compression and the actions of the matrix exponentials. As these quantities are obviously not independent, however, a proper implementation requires some care. We have not used this feature in our numerical experiments and instead rely on experience to choose reasonable tolerances.
\end{remark}

Finally, we present the full procedure for approximating the solution to Equation~\eqref{eq:DRE} in algorithmic form in Algorithms~\ref{alg:expG}--\ref{alg:full}. We consider only the symmetric case of the additive splitting schemes, since the asymmetric version is analogous; change the order of the error estimator from $2s - 2$ to $s-1$ and only use the $L_{+}$ or $L_{-}$ terms instead of both.
When a step is rejected, it is likely that it is because the approximation of $I_Q(h_n)$ is poor. We thus first recompute the whole sequence of quadrature nodes and then retry the step with the same step size. Only if this also fails do we decrease the step size and proceed as normal.

\begin{algorithm}
  \caption{Computing the low-rank factorization of $\cTG(h)P_0$}
  \label{alg:expG}
  \begin{algorithmic}[1]
    \REQUIRE Matrices $S \in \R^{N \times N}$, $L_0 \in \R^{N \times r}$ and $D_0 \in \R^{r \times r}$ with $P_0 = L_0 D_0 L_0^T$, step size $h$
    \STATE Compute $D = (I + hD_0L_0^TSL_0)^{-1}D_0$
    \STATE Set $L = L_0$
    \ENSURE Matrices $L$ and $D$ such that $LDL^T \approx \cTG(h)P_0$
  \end{algorithmic}
\end{algorithm}

\begin{algorithm}
  \caption{Computing the low-rank factorization of $\cTF(h)P_0$}
  \label{alg:expF}
  \begin{algorithmic}[1]
    \REQUIRE Matrices $L_0 \in \R^{N \times r}$, $D_0 \in \R^{r \times r}$ such that $P_0 = L_0 D_0 L_0^T$, step size $h$, approximate low-rank factorization $L_I D_I L_I^T$ of $I_Q(h)$ 
    \STATE Compute $\hat{L} = \exp{hA^T}L_0$
    \STATE Form $L = [\hat{L}, L_I]$ and $D = \begin{bmatrix} D_0 & 0 \\ 0  & D_I\end{bmatrix}$ and column-compress
    \ENSURE Matrices $L$ and $D$ such that $LDL^T \approx \cTF(h)P_0$
  \end{algorithmic}
\end{algorithm}

\begin{algorithm}
  \caption{Approximating the solution to Equation~\eqref{eq:DRE}}
  \label{alg:full}
  \begin{algorithmic}[1]
    \REQUIRE Matrices $A, S \in \R^{N \times N}$, $L_Q \in \R^{N \times r_Q}$ and $D_Q \in \R^{r_Q \times r_Q}$ such that $Q = L_QD_QL_Q^T$,  $L_0 \in \R^{N \times r}$ and $D_0 \in \R^{r \times r}$ such that $P_0 = L_0D_0L_0^T$
    \REQUIRE Desired method order $2s$, coefficients $\{\gamma_k\}_{k=1}^{s}$ for order $2s$, coefficients $\{\beta_k\}_{k=1}^{s}$ for order $2s-2$, initial time step $h_1$, desired error tolerance $\text{TOL}$
    \REQUIRE Equidistant nodes $s_k$ and weights $w_k$ for a quadrature rule of order $s+1$ on $[0, h_1]$ 

    \STATE  Set $\alpha_k = \gamma_k - \beta_k$ for $k = 1, \ldots, s-1$ and $\alpha_s = \gamma_s$
    \STATE Set $k_I =0.2/(2s-2)$, $k_P = 0.2/(2s-2)$
    \STATE Set $n = 1$, $t_n = 0$ and $e_n = 0$

    \WHILE{ $t_n + h_n \le T$ }

    \STATE Low-rank approximate $I_Q(h_n) \approx L_I D_I L_I^T$ according to Algorithm~\ref{alg:update_IQh}, store the computed $L(s_k)$

    \STATE Compute in parallel $L^j_{\pm}$ and $D^j_{\pm}$ such that
    \begin{align*}
      L^j_{+}D^j_{+}(L^j_{+})^T &= \Big(\cTF(h/j)\cTG(h/j)\Big)^j L_{n-1} D_{n-1} L_{n-1}^T \quad \text{and } \\
      L^j_{-}D^j_{-}(L^j_{-})^T &= \Big(\cTG(h/j)\cTF(h/j)\Big)^j L_{n-1} D_{n-1} L_{n-1}^T,
    \end{align*}
    for $j = 1, \ldots, s$, according to Algorithms~\ref{alg:expG} and~\ref{alg:expF}.

    \STATE Form $L_n = \begin{bmatrix}  L_{+}^1 & L_{-}^1 & \cdots & L_{+}^s & L_{-}^s  \end{bmatrix}$, $D_n = \blkdiag(\gamma_1 D^1_{+}, \gamma_1 D^1_{-}, \ldots,  \gamma_s D^s_{+}, \gamma_s D^s_{-})$ \\and column compress 

    \STATE Form  $\hat{L}_n = L_n$, ${\hat{D}_n = \blkdiag(\alpha_1 D^1_{+}, \alpha_1 D^1_{-}, \ldots, \alpha_s D^s_{+}, \alpha_s D^s_{-})}$ and column compress 

    \STATE Compute the local error estimate $e_{n+1} = \norm{\hat{L}_n \hat{D}_n \hat{L}_n^T}_F = \big( \trace \big((\hat{L}_n^T \hat{L}_n \hat{D}_n)^2 \big) \big)^{1/2}$

    \IF{$e_{n+1} > \text{TOL}$}
    \STATE Reject the step 
    \STATE If first rejection, do a full recomputation of $I_Q(h_n)$ and redo step with same $h_n$
    \STATE If still rejected, redo step with $h_n = \Big(\frac{0.9 \text{TOL}}{e_{n+1}} \Big)^{1/(2s-2)} h_n$
    \ELSE
    \STATE Set $t_n = t_{n-1} + h_n$
    \IF{$t_n = T$}
    \BREAK
    \ENDIF
    \STATE     Update the time step by $h_{n+1} = \Big( \frac{ 0.9 \, \text{TOL}}{e_{n+1}}\Big)^{k_I} \Big( \frac{e_n}{e_{n+1}}\Big)^{k_P} h_n$  

    \STATE Set $n = n + 1$
    \ENDIF   

    \IF{$t_n + h_n > T$}
    \STATE Set $h_n = T - t_n$
    \ENDIF

    \ENDWHILE

    \ENSURE Time steps $t_k \in [0, T]$, approximations $L_k$, $D_k$ such that $L_k D_k L_k^T \approx P(t_k)$
  \end{algorithmic}
\end{algorithm}

\section{Numerical experiments} \label{sec:experiments}
In order to verify the validity of the proposed splitting schemes, a number of numerical experiments were performed using MATLAB implementations of the presented algorithms.

Different norms may be used to measure the errors. In all our experiments, we consider relative errors at the final time, measured in the Frobenius norm. That is, if the approximation $P_n$ and a given reference approximation $P_{\text{ref}}$ both approximate the solution $P(T)$, the error is given by
\begin{equation*}
  \frac{\norm{P_n - P_{\text{ref}}}_F}{\norm{P_{\text{ref}}}_F},
\end{equation*}
where $\norm{\cdot}_F$ denotes the Frobenius norm.

\subsection{Order investigation, small-scale}\label{subsec:order_small}
As a first test, we demonstrate that the methods exhibit the expected orders of convergence when constant step sizes are used.  For this, we consider a small-scale problem with $N = 10$ and take $A$, $Q$, $S$ and $P_0$ to be random matrices with the latter three having rank $4$. The small dimension of the problem means that we may compute a highly accurate reference approximation by unrolling the matrix-valued problem into a vector-valued problem of dimension $N^2$ and applying a standard method for ODEs. Here we utilize the MATLAB built-in function \texttt{ode15s}, which implements an adaptive variable-order multistep method, with an absolute tolerance of $10^{-20}$ and a relative tolerance of ${2.22 \cdot 10^{-14}}$ (the minimum).

For this test, we consider the asymmetric splitting schemes~\eqref{eq:asym} of orders $2$ and $3$, the symmetric schemes~\eqref{eq:sym} of orders $2$, $4$, $6$ and $8$, as well as the 2nd-order Strang splitting. To compute terms of the form $\exp{hA^T}L$, we use the 5th-order implicit Runge-Kutta scheme RadauIA~\cite[Chapter IV.5]{HairerWannerII} and halve the step size until two subsequent approximations differ (relatively) by at most $10^{-6}$. This can clearly be done better, ideally with adaptive time stepping also on this level, but it is sufficient for our purposes. We set the column compression tolerance to $10^{-16}$ so that it has no effect on the results.

\begin{figure}[!t] 
  \centering{ \includegraphics[width=1.0\columnwidth]{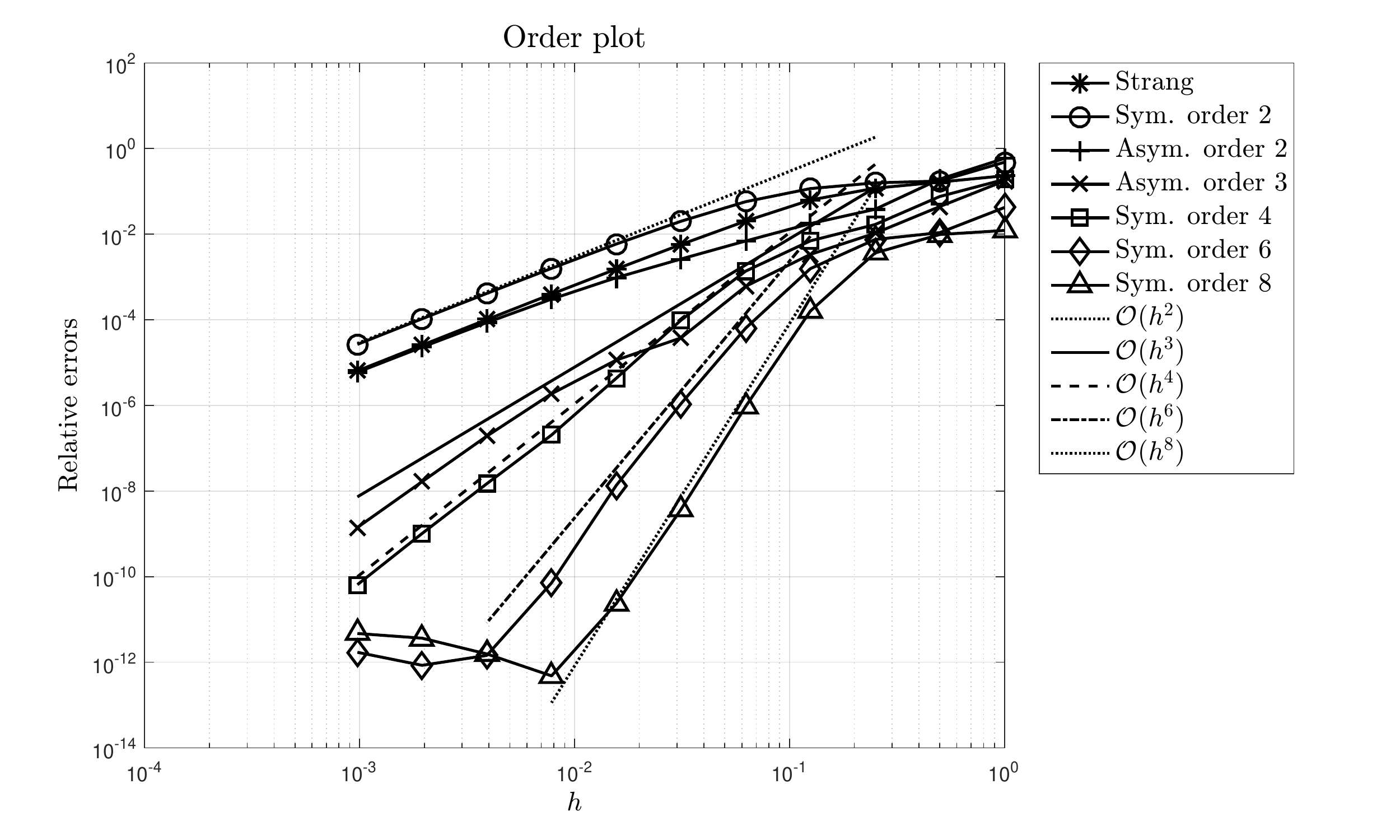} }
\caption{Errors plotted against step sizes for the problem defined in Section~\ref{subsec:order_small}. We observe that all the methods exhibit the expected convergence orders until the round-off level is reached, except for very large step sizes. }
\label{fig:order_small} 
\end{figure}

The results are shown in Figure~\ref{fig:order_small}, where it can be seen that all the methods do, indeed, achieve the expected converge orders. However, a few comments are in order. First, the 3rd-order asymmetric scheme actually exhibits an order of convergence which is slightly larger than $3$. This is not true in general and we interpret this as the structure of the error being favourable for this particular problem. Secondly, the errors for the 6th- and 8th-order methods level out around $10^{-12}$. This is due to round-off error accumulation in each step. Using a dense instead of low-rank factored version of the code, computing $\exp{hA}$ explicitly and approximating $I_Q(h)$ to high accuracy gives similar results. The leveling out of all the error curves for large step sizes is due to leaving the asymptotic regime; for these step sizes also lower-order error terms influence the result. Thirdly, we note that the 2nd-order asymmetric method performs slightly better than both the Strang splitting and the 2nd-order symmetric method. However, since it is $50 \%$ more expensive if parallelization is not used, and even more so if it is, we clearly still prefer the symmetric method. 

\subsection{Order investigation, larger-scale }\label{subsec:order_large}
We consider also a larger, real-world problem, arising from the optimal control of steel cooling~\cite{BennerSaak2005,Saak2003}. This is essentially a finite-element discretization of a semi-linear PDE given on a non-convex two-dimensional domain. It results in matrices $A \in \R^{N \times N}$, $B \in \R^{N \times 7}$ and $C \in \R^{6 \times N}$ from which we construct $Q = C^TC$ and $S = B R^{-1} B^T$, with $R^{-1} = I$. The problem also involves a mass matrix, i.e.\ the state equation is $M\dot{x} = Ax + Bu$. We handle this without inverting $M$ by straightforward modifications to the code as in~\cite{Stillfjord2015}. Additionally, due to a scaling of the problem, a simulation time step of $1$ second corresponds to a real time step of $10^{-2}$ seconds. To avoid confusion, we work with the simulation time throughout, and therefore use a final time $T = 4500$.

The exact solution to the problem is unavailable, and since the other currently existing methods are limited to low orders it is infeasible to use these to compute a sufficiently accurate reference approximation. Instead, we use the 8th-order symmetric splitting scheme itself for this, but with a step size half as large as the smallest step size for the actual approximations. In this experiment we do employ parallelization through use of MATLAB's \texttt{parfor} command, using 8 cores on a cluster built out of Intel 2650v3 CPUs. We restrict ourselves to the Strang splitting and the symmetric methods, since our tests indicate that these are typically more efficient than the asymmetric methods. We perform two tests, one with $N = 371$ and one with $N = 1357$. Except the time step size, the only varying parameter is the relative tolerance for computing the matrix exponential actions. In the smaller example, this is set to $10^{-3}$ for the Strang splitting and $10^{-3}$, $10^{-6}$, $10^{-8}$ and $10^{-8}$ for the additive schemes of order $2$, $4$, $6$ and $8$, respectively. In the larger example, we take instead $10^{-3}$, $10^{-3}$, $10^{-5}$, $10^{-6}$ and $10^{-6}$, respectively. The column compression tolerance is in all cases set to $N\epsilon$, where $\epsilon$ is the machine epsilon.

\begin{figure}[!t] 
    \centering{
      \includegraphics[width=0.48\columnwidth]{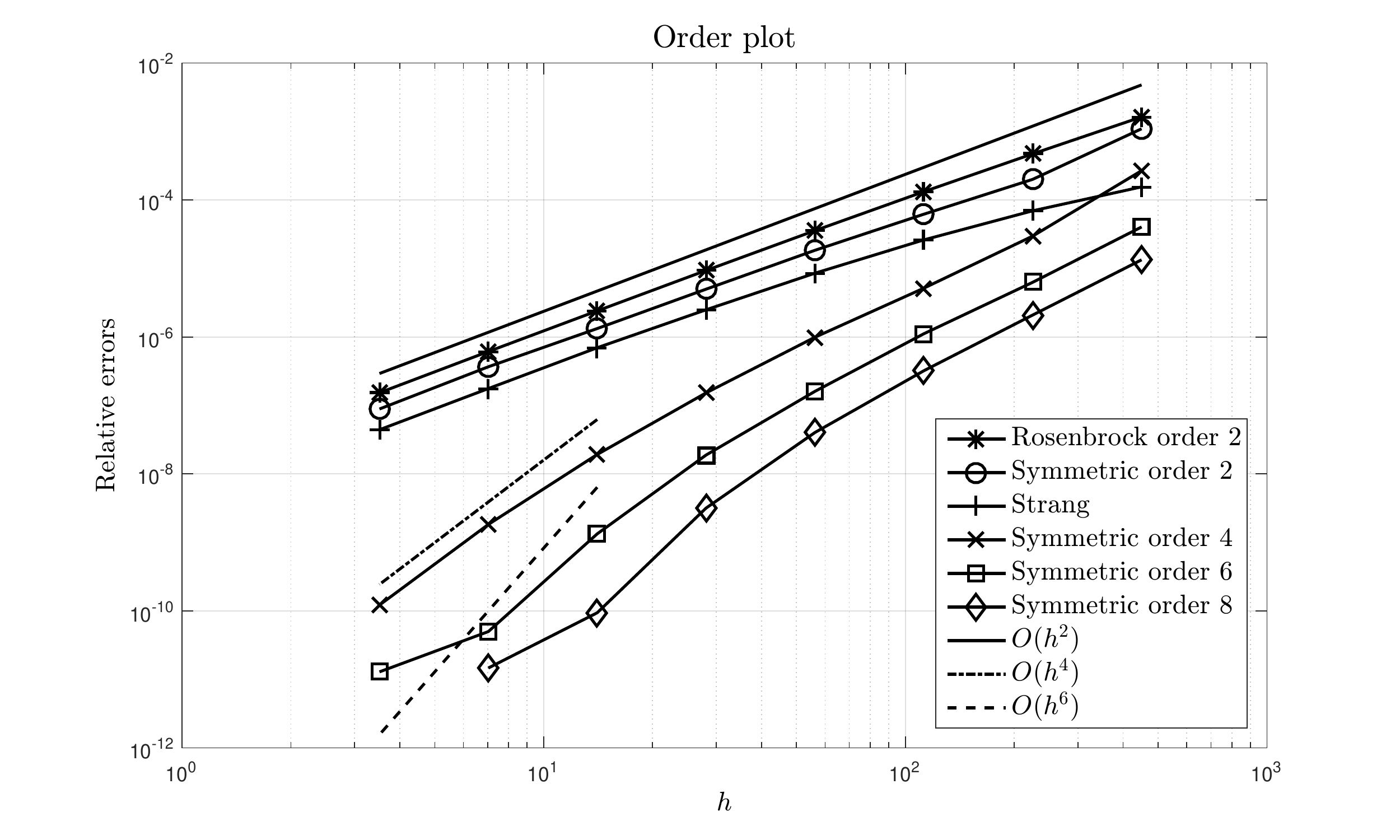}
      \includegraphics[width=0.48\columnwidth]{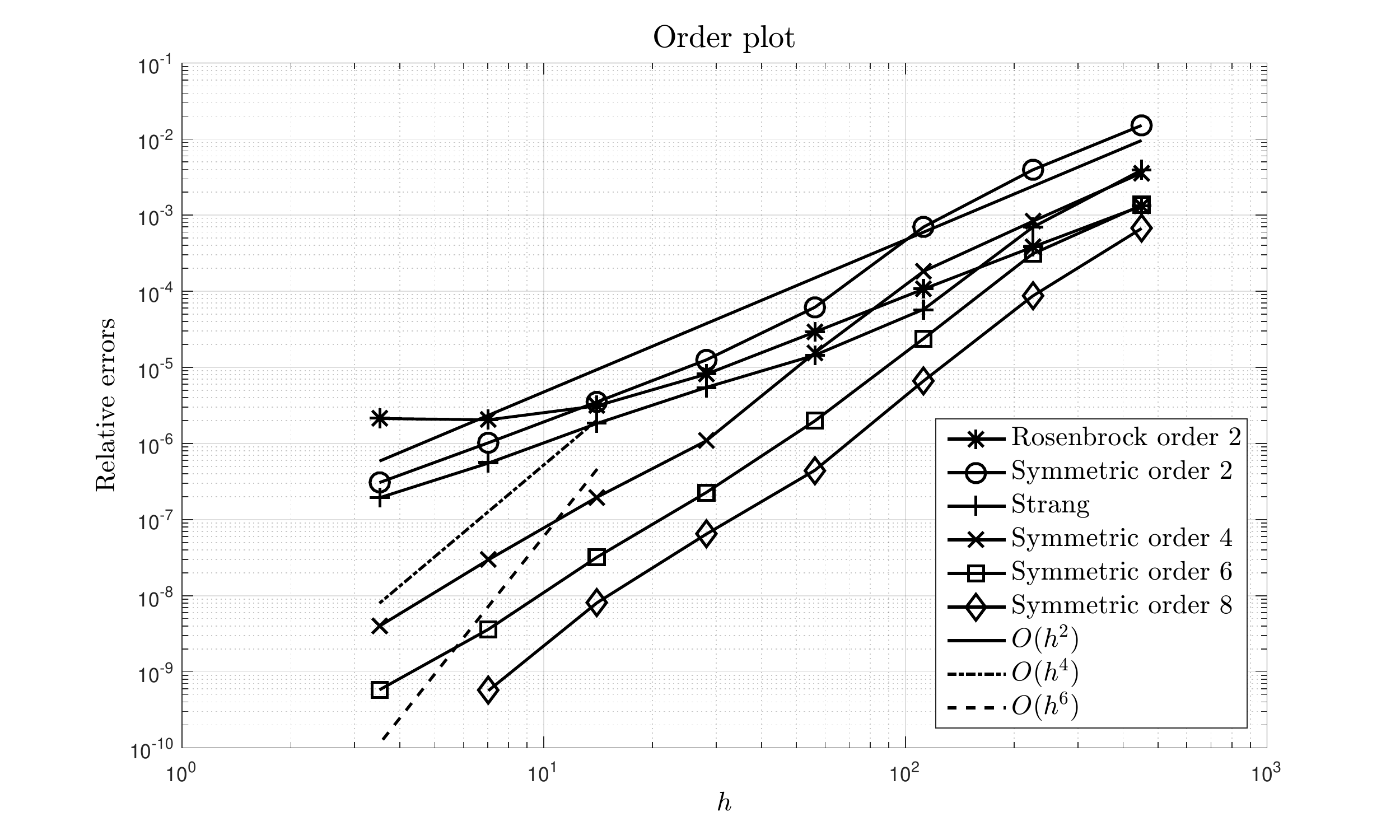} 
    }
\caption{Errors plotted against step sizes for the problem defined in Section~\ref{subsec:order_large}. Left: $N = 371$. Right: $N = 1357$. We observe that the second-order methods show second-order behaviour for all step sizes used, while the higher-order methods suffer from order reduction. For the smaller problem size, we recapture the higher-order behaviour for the smallest step sizes, while the larger problem size requires even smaller step sizes before this happens. Regardless of this, the errors of the higher-order methods are significantly smaller than those of the second-order methods.}
\label{fig:order_large} 
\end{figure}

 Figure~\ref{fig:order_large} shows the results, with $N = 371$ on the left and $N = 1357$ on the right. In the smaller example, we observe that the second-order methods behave as expected, while the higher-order methods only achieve their respective orders for small step sizes. The fact that the errors level out at around $10^{-11}$ can be avoided by computing the matrix exponentials more accurately, but at additional cost.
In the larger example, the situation is slightly worse in that neither of the higher-order methods reach their asymptotic regimes with the used step sizes. This issue may be due to a lack of regularity in the solution to the exact problem. As in the smaller example, we could eliminate the leveling out of the error by decreasing the tolerance for the matrix exponential actions. However, as the computation times required for these small errors are already rather long, we do not do this. In spite of these issues, we note that the higher-order methods still produce much smaller errors for all step sizes except the largest.

Also included in Figure~\ref{fig:order_large} are the corresponding errors for the second-order Rosenbrock method proposed in~\cite{BennerMena2013}. These computations were done using the MATLAB software M-MESS 1.0.1~\cite{MMESS101}, which implements the improved $LDL^T$-formulation given in~\cite{LangMenaSaak2015}. We choose the parameters suggested in the example code for the steel cooling problem. For the smaller problem, we observe clear second-order convergence with errors of comparable size to the splitting schemes. The situation is similar in the larger problem, except that the error evens out for small step sizes, likely due to inner iterations not being computed accurately enough.

\begin{figure}[!t] 
    \centering{
      \includegraphics[width=0.48\columnwidth]{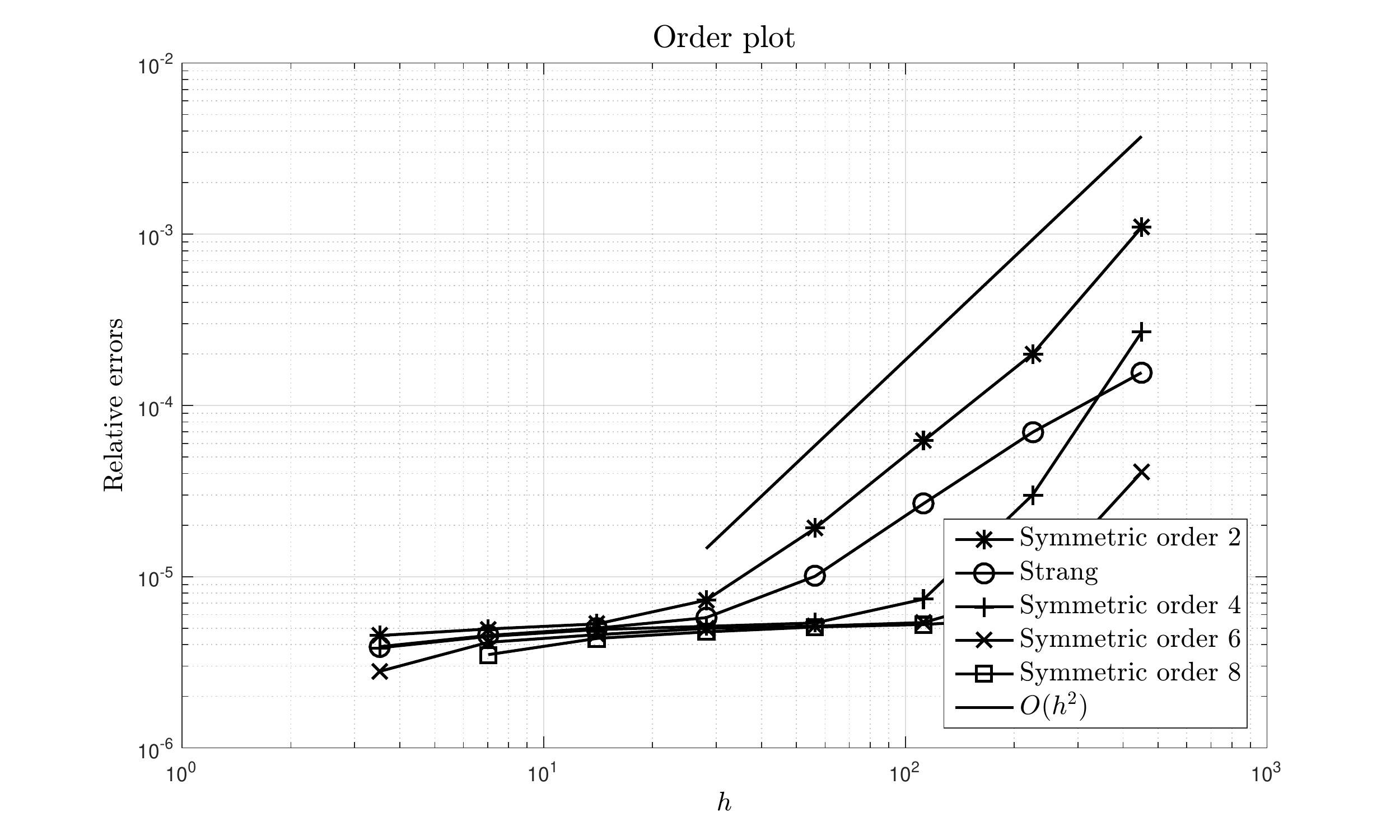}
      \includegraphics[width=0.48\columnwidth]{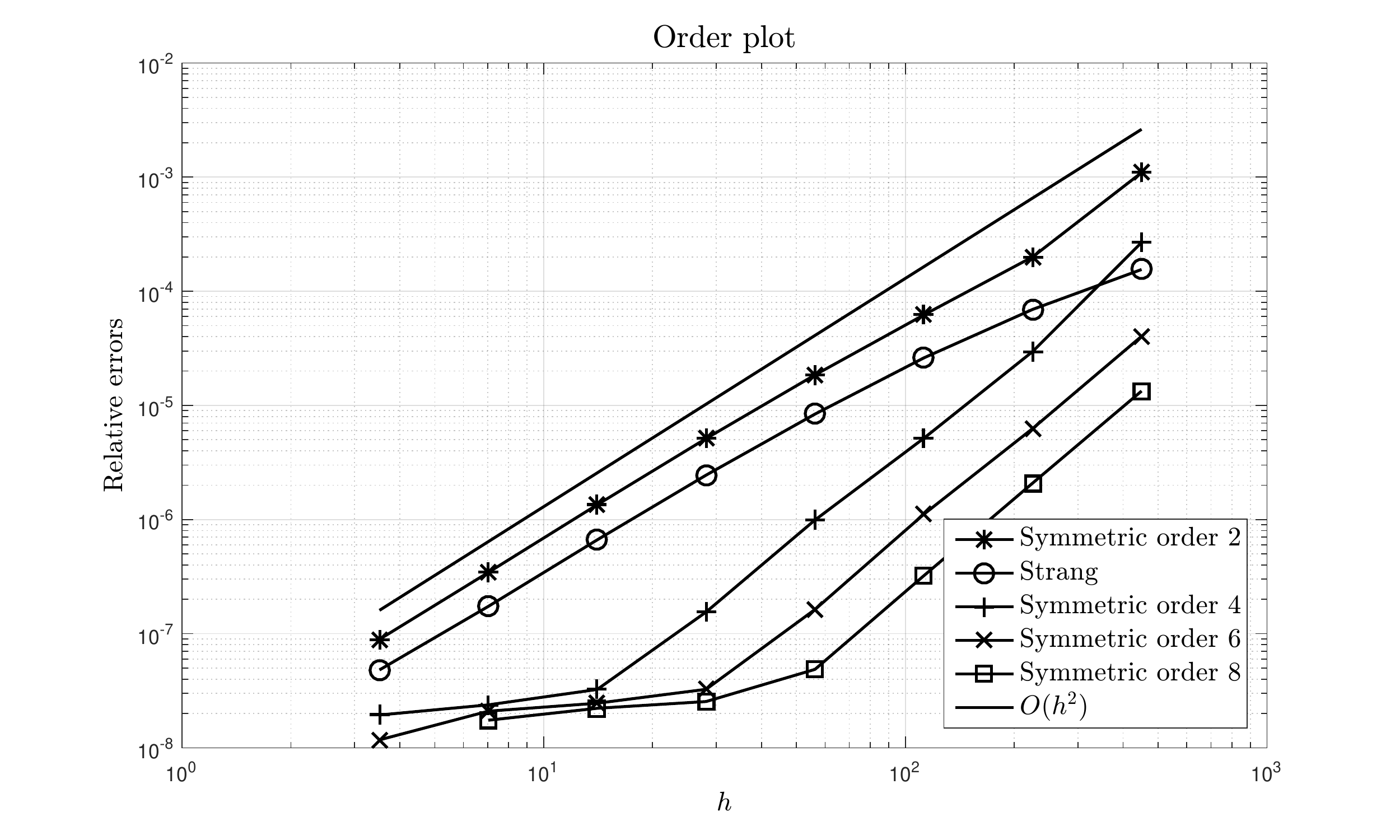} 
    }
\caption{Errors plotted against step sizes for the problem defined in Section~\ref{subsec:order_large} with $N = 371$ and the column compression tolerances $10^{-8}$ (left) and $10^{-10}$ (right). We note that the convergence is unaffected until the truncation level is reached. Because these errors are introduced in each time step, the error levels out at a value larger than the specified tolerance.}
\label{fig:compression} 
\end{figure}

Finally, we note that the column compression tolerance has been chosen rather small. This is required for the small step sizes, due to the small errors produced by the high-order methods. In Figure~\ref{fig:compression}, we demonstrate the effect of increasing this tolerance. We note that the convergence behaviour is unaffected until the truncation level is reached.
Obviously, the ranks of the approximations are heavily affected by changing this tolerance. Table~\ref{table:ranks} illustrates this, by tabulating the rank of the Strang splitting approximation at the final time when $N=371$. The ranks of the other methods differ (at most) by $\pm 10$ from these values for the two largest step sizes, and by $\pm 3$ for the other step sizes. In all cases, the rank increases monotonically until the final time, i.e.\ the presented ranks are the maximum attained during the simulation. For comparison, the rank of the corresponding ARE (which the DRE solution tends to as $t \to \infty$) is $138$. When $N=1357$, the ranks of the approximations are of similar size, which fits well with the expectation that low rank is a property inherent to the DRE, independent of the discretizations.

\begin{table}                                                 
\centering                                                    
\begin{tabular}{c|cccccccc}                                    
  Tolerance \textbackslash \; No.\ steps  &  10 &  20 &  40 &  80 & 160 & 320 & 640 & 1280 \\
  \hline
  \rule{0pt}{2.6ex} 
  $10^{-8}$           &  68 &  70 &  71 &  70 &  68 &  67 &  66 &   65 \\
  $10^{-10}$          &  83 &  86 &  86 &  86 &  85 &  84 &  84 &   82 \\
  $8.2\cdot 10^{-14}$ & 102 & 107 & 109 & 110 & 110 & 109 & 107 &  107 \\  
\end{tabular}                                                 
\caption{The ranks of the Strang splitting approximation at the final time when $N = 371$, for different column compression tolerances and different number of time steps. The last value $8.2\cdot 10^{-14}$ is equal to $N\epsilon$, where $\epsilon$ is the machine epsilon.}
\label{table:ranks}                                    
\end{table}

\subsection{Efficiency} \label{subsec:efficiency}
\begin{figure}[!t] 
  \centering{ \includegraphics[width=1.0\columnwidth]{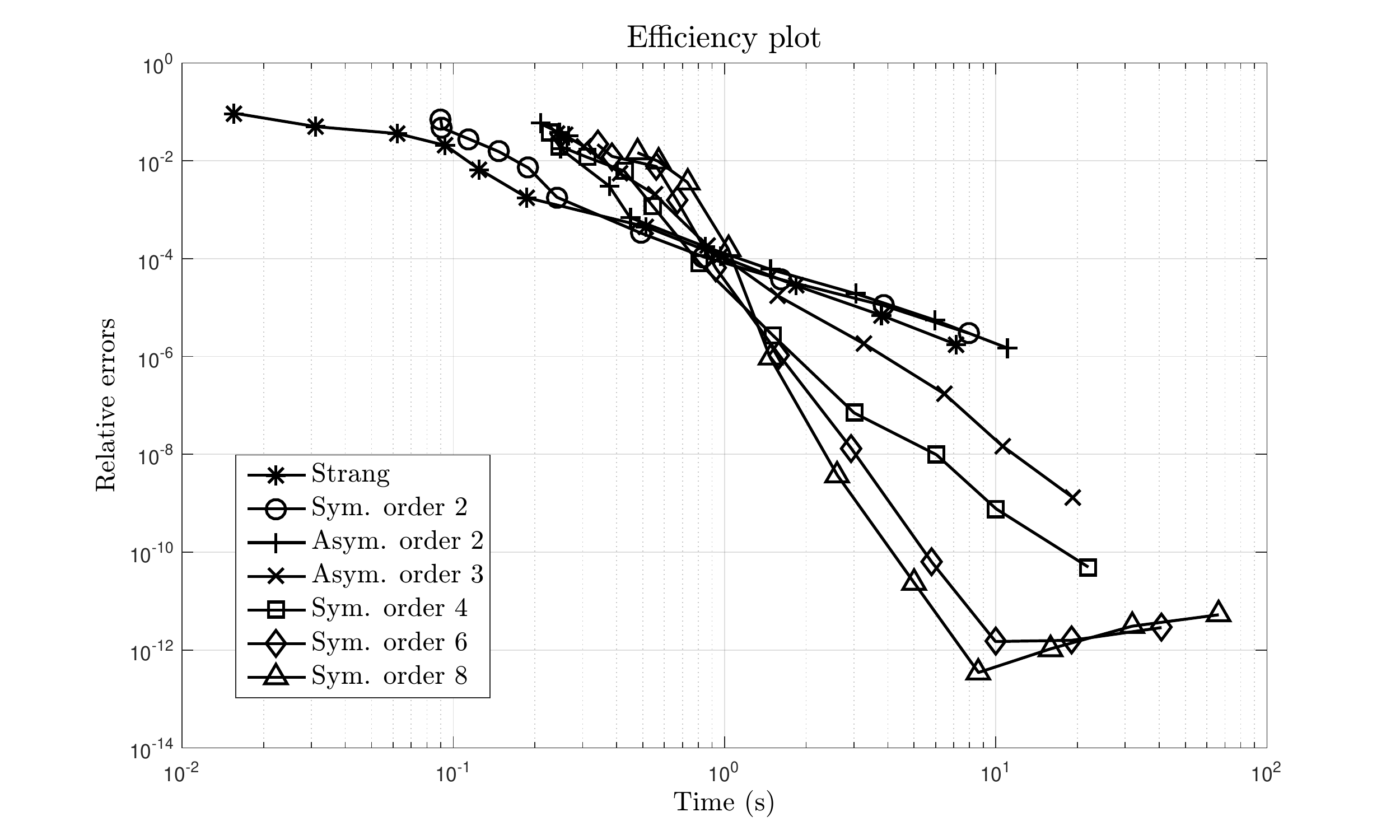} }
\caption{Errors plotted against computation times for the problem defined in Section~\ref{subsec:order_small}. We see that the lower-order methods are most efficient for high error levels, while the higher-order methods are most efficient for low error levels. For errors around $10^{-4}$, the efficiency of all the methods is comparable.}
\label{fig:eff_small} 
\end{figure}

While the higher-order methods produce smaller errors, this is only relevant if their computational costs are similar to that of the lower-order methods.
We therefore also provide a rough comparison of the efficiency of the different methods. Figure~\ref{fig:eff_small} shows the errors plotted against the required computation time (wall-clock time) for the small-scale problem given in Section~\ref{subsec:order_small}, with the same method parameters. These are the same errors as in Figure~\ref{fig:order_small}, i.e.\ the step size $h$ is the only varying parameter. We observe that all the methods are roughly equivalent for high tolerances, while for error levels below $10^{-4}$, the symmetric methods outperform the others. For very small errors, the 6th- and 8th-order methods are clearly superior. This is in spite of the fact that parallelization was not used in this case (since the extra time spent on transferring data was much larger than the actual computation time).

\begin{figure}[!t] 
  \centering{ \includegraphics[width=0.48\columnwidth]{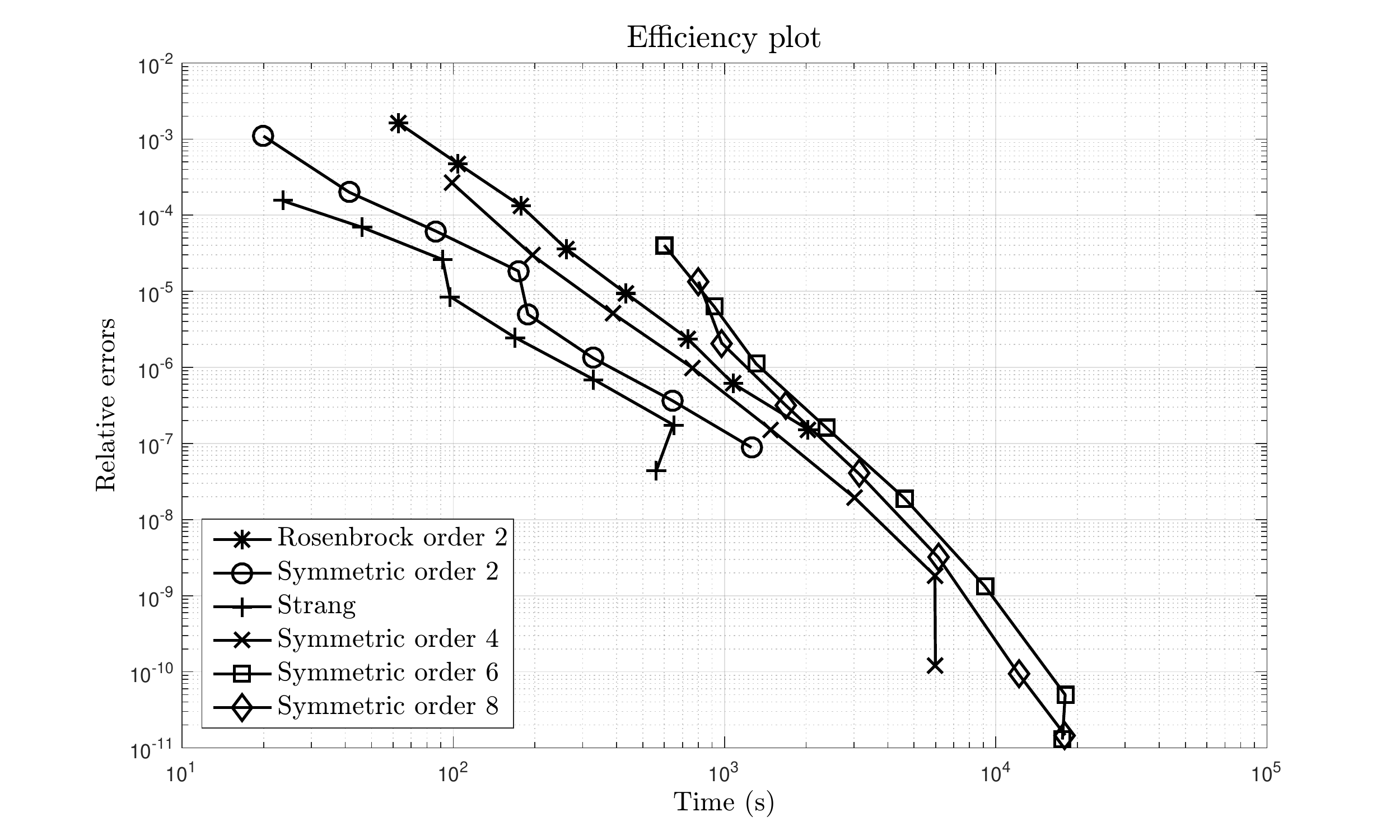}
    \includegraphics[width=0.48\columnwidth]{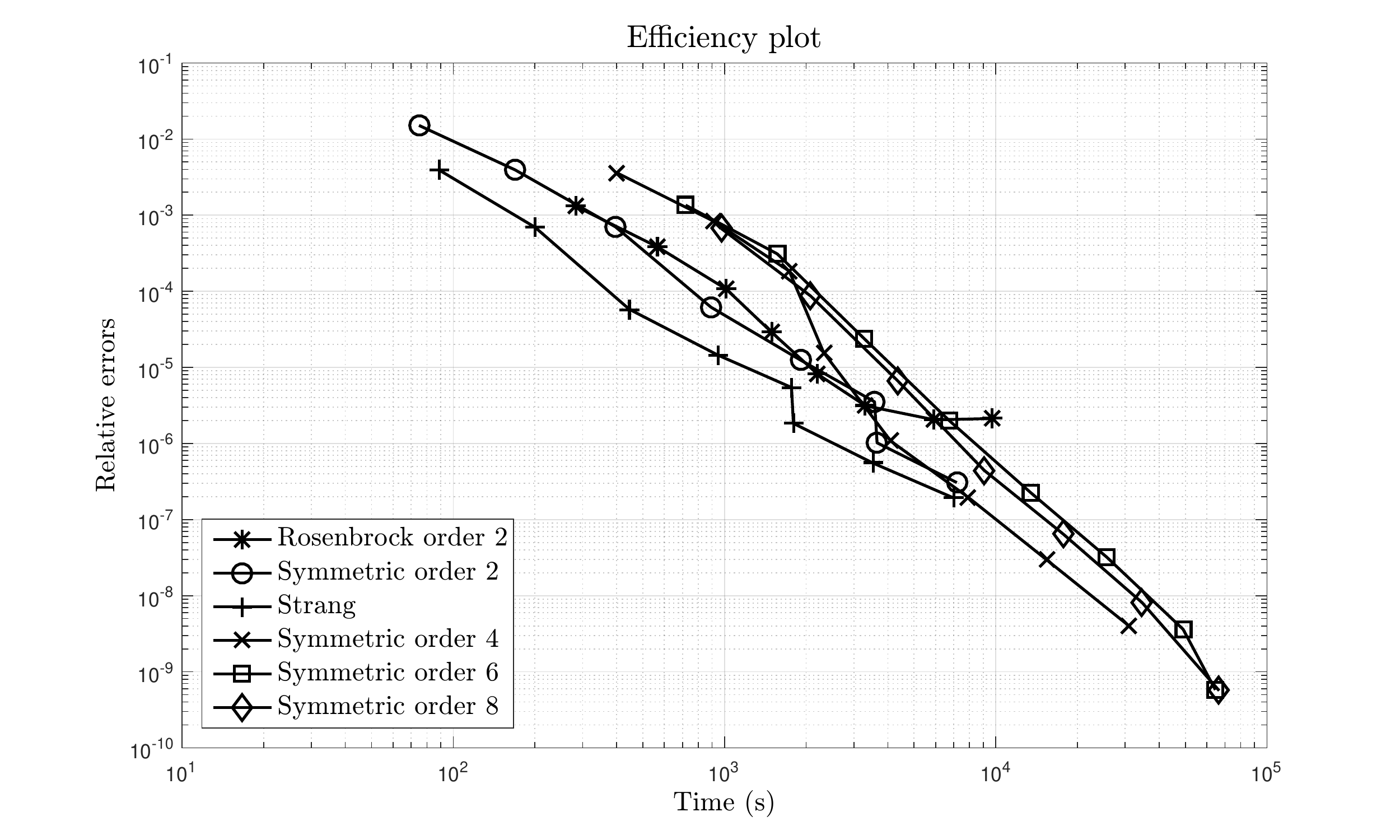} }
\caption{Errors plotted against computation times for the problem defined in Section~\ref{subsec:order_large}. Left: $N = 371$. Right: $N = 1357$. The lower-order methods are again most efficient for high error levels and vice versa, though the difference between the methods is much less than in Figure~\ref{fig:eff_small}.}
\label{fig:eff_large} 
\end{figure}

The results for the steel cooling problem are shown in Figure~\ref{fig:eff_large}. These are similar to the small-scale case in that the higher-order methods are more efficient for small errors while the Strang splitting is most efficient for large errors. The plot is slightly misleading, because the low matrix exponential tolerances required for the high-order methods to reach the smallest errors are not strictly required for the less accurate approximations. Similarly, the lower-order methods would need to compute the matrix exponentials more accurately when the step size is further decreased. Thus the real cut-off point where the higher-order methods become more efficient lies somewhere between the error levels $10^{-5}$ and $10^{-7}$. We also observe that the 8th-order method is superior to the 6th-order method for small errors. This is due to the parallelization: the cost of increasing the order by 2 is equivalent to only one extra Lie splitting step, and one extra processor. Using even higher orders may thus be beneficial, but eventually the overhead costs incurred by the parallelization will dominate.

 The strange kinks in the error curves require an explanation. For the Strang splitting this happens twice when $N = 371$, and on the latter occurrence the computation time even decreases slightly when the step size is decreased. This happens due to the way we compute the actions of the matrix exponentials: if the requested accuracy is not reached, the computation is repeated with twice as many sub-steps. Reducing the time step by a factor two makes this computation easier, and it may thus be that most of these computations require only half as many sub-steps as for the larger time step. With twice as many time steps, the total computation time is therefore roughly unchanged.

 To provide an indication of what parts of the splitting schemes are expensive, all the methods were run through MATLAB's \texttt{profile} command while solving the steel cooling problem with $N=1357$ and with either $40$ or $640$ time steps, corresponding to $h = 112.5$ or $h = 7.0313$. The tolerance for the $\exp{\gamma h A^T}L$-computations was set to $10^{-6}$ in all cases. The results are shown in Table~\ref{table:breakdown_40} and~~\ref{table:breakdown_640}. We observe that, as expected, the evaluations of $\cTG$ are essentially free in comparison to $\cTF$. The cost of the latter completely dominates the overall procedure. Along with the observation in the previous paragraph, this provides additional incentive for studying better implementation strategies for this basic operation.

Further, we observe that the relative cost of column compression increases with the order of the method, since $L^j_{\pm}$ and $D^j_{\pm}$ increase in size. In total, however, this cost is negligible, despite the fact that the ranks of the approximations are not very small. It should be noted here that due to the difficulties of accurately timing parallel code in this level of detail, a serial implementation was used. While this skews the ratios, the effect is very small because almost all column compressions originate from $\cTF$-evaluations ($95$\% for the order $8$ method and the smallest step size). 
As expected, the relative cost of the one-time computation of $I_Q(h)$ is higher for a small number of time steps, but even in the worst case it is measured in single digit percentages. With many time steps, the relative cost is negligible.

\begin{table}                                                   
\centering                                                      
\begin{tabular}{c|ccccc}                                         
Operation \textbackslash \; Method & Strang & Additive 2 & Additive 4 & Additive 6 & Additive 8 \\
\hline
$\cTG$ & 0.01 & 0.03 & 0.03 & 0.03 & 0.03 \\
$\cTF$ & 96.60 & 95.90 & 97.50 & 98.06 & 98.32 \\               
Column compression & 0.16 & 0.31 & 0.33 & 0.36 & 0.38 \\        
$\exp{\gamma h A^T}L$ & 99.10 & 98.88 & 99.16 & 99.16 & 99.18 \\
$I_Q(h)$ & 3.38 & 3.92 & 2.30 & 1.71 & 1.42 \\                  
\end{tabular}                                                   
\caption{Computational time breakdown for the splitting schemes when applied to the steel cooling problem with $N = 1357$ and $40$ time steps. Shown is the time spent on the given operation, divided by the total time for the integration (in percent). The numbers are not independent, e.g.\ computing $I_Q(h)$ requires several $\exp{\gamma h A^T}L$ evaluations and column compressions.}
\label{table:breakdown_40}                                      
\end{table}

\begin{table}                                                
\centering 
\begin{tabular}{c|ccccc}                                         
Operation \textbackslash \; Method & Strang & Additive 2 & Additive 4 & Additive 6 & Additive 8 \\
\hline
$\cTG$ & 0.02 & 0.03 & 0.03 & 0.03 & 0.03 \\                    
$\cTF$ & 99.86 & 99.66 & 99.68 & 99.68 & 99.67 \\               
Column compression & 0.13 & 0.31 & 0.33 & 0.34 & 0.37 \\        
$\exp{\gamma h A^T}L$ & 99.51 & 99.33 & 99.30 & 99.30 & 99.27 \\ 
$I_Q(h)$ & 0.12 & 0.13 & 0.08 & 0.05 & 0.04 \\                  
\end{tabular}            
\caption{Computational time breakdown for the splitting schemes when applied to the steel cooling problem with $N = 1357$ and $640$ time steps. Shown is the time spent on the given operation, divided by the total time for the integration (in percent). The numbers are not independent, e.g.\ computing $I_Q(h)$ requires several $\exp{\gamma h A^T}L$ evaluations and column compressions.}
\label{table:breakdown_640}
\end{table}

Finally, we note that Figure~\ref{fig:eff_large}, like Figure~\ref{fig:order_large}, also includes the results for the second-order Rosenbrock method. While a fair comparison is difficult, and many parameters could be further fine-tuned for all the methods, these results clearly indicate that the splitting schemes constitute a competitive alternative to this class of methods.

\subsection{Time adaptivity} \label{subsec:adaptivity}
\begin{figure}[!t] 
  \centering{ \includegraphics[width=1.0\columnwidth]{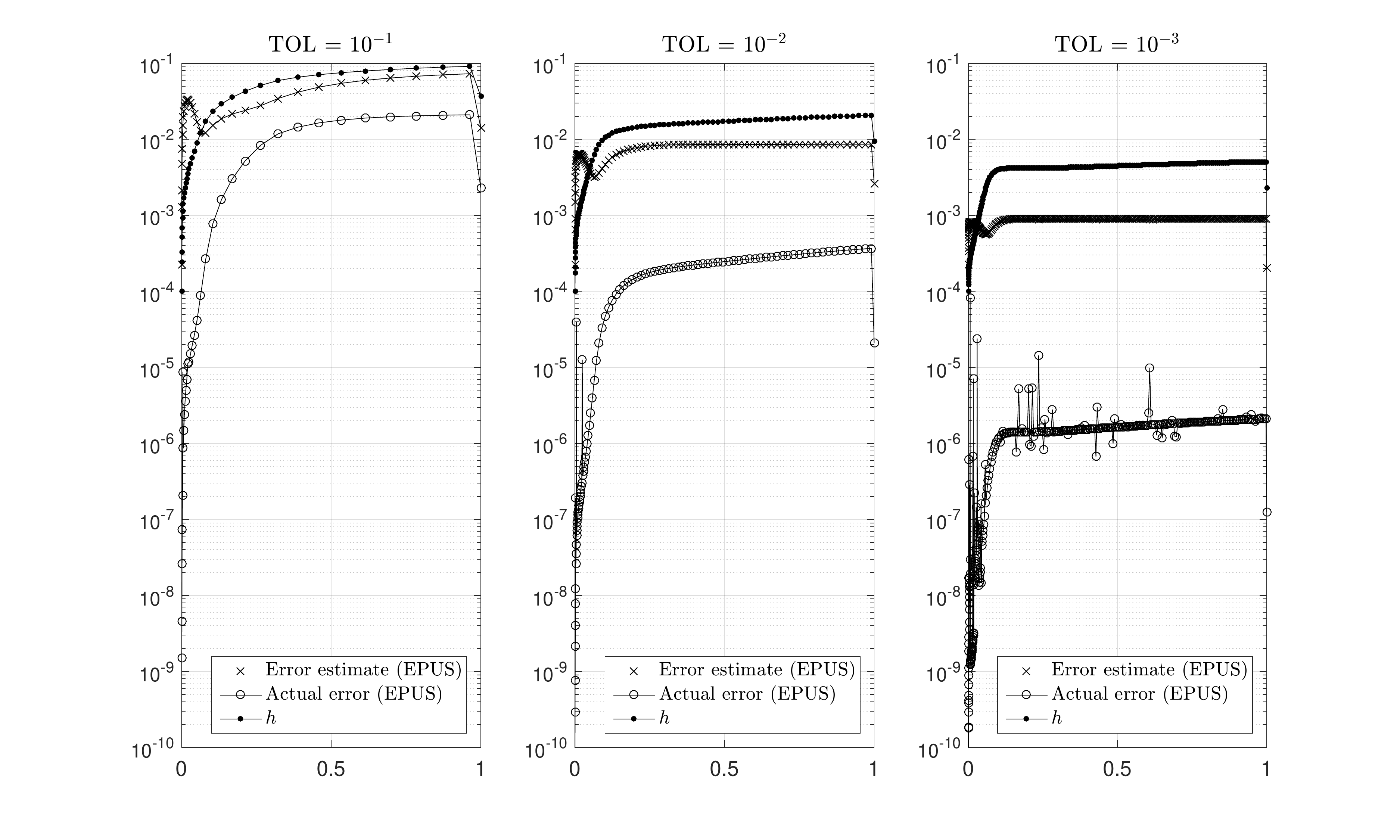}}
\caption{Computed error estimate, actual error and step size for each time step, when applying the adaptive 4th-order symmetric splitting scheme to the problem defined in Section~\ref{subsec:order_small}. We consider error per unit step (EPUS), i.e. all errors are divided by the time step. The tolerances used are, from left to right, $10^{-1}$, $10^{-2}$ and $10^{-3}$. We observe that the adaptivity finds the maximum step size such that the error estimate is equal to the tolerance. Due to the lower-order estimate, the actual error is in all cases less than the estimated error, and the difference increases as the step size decreases. The sudden drop in step size (and error) in the final step is necessary in order to exactly reach the final time $T$.}
\label{fig:adaptive_small} 
\end{figure}

Finally, we test the full time step adaptive code with the 4th-order symmetric splitting scheme. In Figure~\ref{fig:adaptive_small} we have plotted the results of using four different tolerances on the small-scale problem defined in Section~\ref{subsec:order_small}. We plot both the error estimated by the method using the embedded method, and the actual error. The latter is computed by using the same method, but by taking $10$ equidistant steps in each of the steps given by the adaptive code. We observe that the actual error is in all cases less than the estimated error, and the difference increases as the tolerance decreases. This is due to the fact that the error estimate is of a lower order than the actual method used. The effect is more pronounced here than usual, since in the symmetric case the accuracy of the estimate is $2$ orders less than the method. In each figure we have also plotted the step sizes, and we see that the controller works well in finding the maximum possible step size. For the largest tolerance, the controller is too cautious and does not quite reach the tolerance until the simulation is over. Effects like this can (and should, this is one area we aim to pursue in the near future) be tuned by adjusting the parameters $k_I$ and $k_P$.

\begin{figure}[!t] 
  \centering{ \includegraphics[width=1.0\columnwidth]{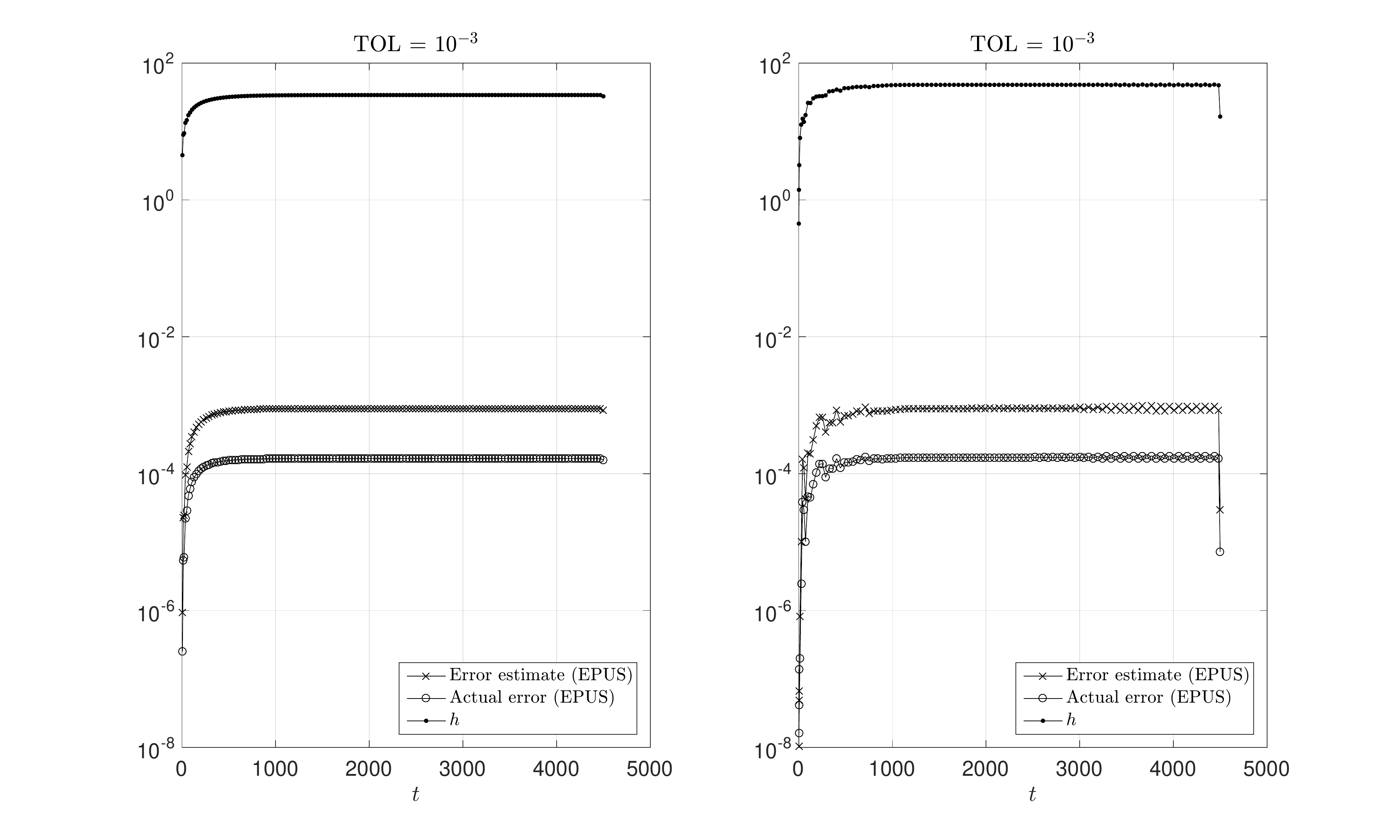}}
\caption{Computed error estimate, actual error and step size for each time step, when applying the adaptive 4th-order symmetric splitting scheme with tolerance $10^{-3}$ to the problem defined in Section~\ref{subsec:order_large}. The method was applied either without (left) or with (right) Algorithm~\ref{alg:update_IQh}. We consider error per unit step (EPUS), i.e. all errors are divided by the time step. We observe that the adaptivity works rather well.}
\label{fig:adaptive_large} 
\end{figure}

In Figure~\ref{fig:adaptive_large}, we have repeated the same experiment but on the steel problem with $N = 371$ and only with the tolerance $10^{-3}$. Also in this case, the adaptiveness seems to work well -- the maximum possible step size (given the tolerance) is quickly reached and after this it varies very little. In this case, the difference between the error estimate and the actual error is not as large as in the previous example. This is likely due to the order reductions observed in Section~\ref{subsec:order_large}.

The left plot shows the results when Algorithm~\ref{alg:update_IQh} is not used and the right plot when it is. In both cases, we used a column compression tolerance of $10^{-8}$, a relative tolerance of $10^{-4}$ for the matrix exponential actions and quadrature of order $9$ to compute $I_Q(h_n)$. When not using Algorithm~\ref{alg:update_IQh} we use Gaussian quadrature rather than Newton-Cotes, and thus these computations only need $4$ quadrature nodes compared to the updating formula which needs $10$. Still, as demonstrated by the computational time breakdown in Table~\ref{table:adaptive_breakdown}, the latter is more efficient because typically only one or even none of these nodes need to be updated in each step. In the current experiment, $111$ steps were taken. Of these, $6$ were rejected which required all $10$ nodes to be updated. During the remaining $105$ steps, a total of $10$ nodes required an update.

\begin{table}                                                      
\centering                                                         
\begin{tabular}{r|ccccc}                                            
Method \textbackslash \; Operation & Total time & $I_Q(h_n)$ & $\cTF$ & $\exp{\gamma h A^T}L$ & CC \\[.4ex]
\hline
  \rule{0pt}{2.6ex} 
Without Algorithm 1 & 100.00 & 46.57 & 52.42 & 87.32 & 1.78 \\     
With Algorithm 1 & 65.90 & 7.97 & 89.90 & 86.90 & 6.14 \\   
\end{tabular}                                                      
\caption{Computational time breakdown for the adaptive splitting scheme when applied to the steel cooling experiment in Section~\ref{subsec:adaptivity}, either with or without the use of Algorithm~\ref{alg:update_IQh}. The first column shows the relative computation times depending on this choice. The other columns show the time spent on the given operation, divided by the total time for the respective method. All numbers are in percent, and ``CC'' is an abbreviation for column compression. Only the numbers in the last two columns are independent. The remaining computation time was spent on (unoptimized) caching of matrices and general bookkeeping. }
\label{table:adaptive_breakdown}                                         
\end{table}

\section{Conclusions} \label{sec:conclusions}
We have introduced a family of splitting schemes for differential Riccati equations which may be of arbitrarily high order, and shown that they may be implemented efficiently in a large-scale setting by utilizing the low-rank $LDL^T$-factorization. Our numerical experiments indicate that the higher-order methods are more efficient when high accuracy is desired, though this of course depends on the actual problem. In addition, we have demonstrated that these methods contain natural embedded error estimates, which e.g.\ may be used for time step adaptivity. While further research on appropriate controller parameters in this setting is required, experiments show that even a basic implementation gives promising results.


\begin{thebibliography}{10}
\providecommand{\url}[1]{{#1}}
\providecommand{\urlprefix}{URL }
\expandafter\ifx\csname urlstyle\endcsname\relax
  \providecommand{\doi}[1]{DOI~\discretionary{}{}{}#1}\else
  \providecommand{\doi}{DOI~\discretionary{}{}{}\begingroup
  \urlstyle{rm}\Url}\fi

\bibitem{AbouKandil_etal2003}
Abou-Kandil, H., Freiling, G., Ionescu, V., Jank, G.: Matrix {R}iccati
  equations.
\newblock Systems \& Control: Foundations \& Applications. Birkh\"auser, Basel
  (2003).
\newblock \doi{10.1007/978-3-0348-8081-7}

\bibitem{AmodeiBuchot2010}
Amodei, L., Buchot, J.M.: An invariant subspace method for large-scale
  algebraic {R}iccati equation.
\newblock Appl.\ Numer.\ Math. \textbf{60}(11), 1067--1082 (2010).
\newblock \doi{10.1016/j.apnum.2009.09.006}

\bibitem{AntoulasSorensenZhou2002}
Antoulas, A.C., Sorensen, D.C., Zhou, Y.: On the decay rate of {H}ankel
  singular values and related issues.
\newblock Syst.\ Control Lett. \textbf{46}(5), 323--342 (2002).
\newblock \doi{10.1016/S0167-6911(02)00147-0}

\bibitem{BasarBernhard1995}
Ba{\c{s}}ar, T., Bernhard, P.: {$H^{\infty}$}-optimal control and related
  minimax design problems, second edn.
\newblock Sys.\ Con.\ Fdn. Birkh\"auser Boston, Inc., Boston, MA (1995).
\newblock \doi{10.1007/978-0-8176-4757-5}.
\newblock A dynamic game approach

\bibitem{BennerBujanovic2016}
Benner, P., Bujanovi{\'c}, Z.: On the solution of large-scale algebraic
  {R}iccati equations by using low-dimensional invariant subspaces.
\newblock Linear Algebra Appl. \textbf{488}, 430--459 (2016).
\newblock \doi{10.1016/j.laa.2015.09.027}

\bibitem{BennerMena2013}
Benner, P., Mena, H.: Rosenbrock methods for solving {R}iccati differential
  equations.
\newblock IEEE Trans.\ Automat.\ Control \textbf{58}(11), 2950--2956 (2013).
\newblock \doi{10.1109/TAC.2013.2258495}

\bibitem{BennerMena2016}
Benner, P., Mena, H.: Numerical solution of the infinite-dimensional {LQR}
  problem and the associated {R}iccati differential equations.
\newblock J.~Numer.\ Math.  (2016).
\newblock \doi{10.1515/jnma-2016-1039}.
\newblock Advance online publication, retrieved 23 Nov. 2016

\bibitem{BennerSaak2005}
Benner, P., Saak, J.: A semi-discretized heat transfer model for optimal
  cooling of steel profiles.
\newblock In: Dimension Reduction of Large-Scale Systems, \emph{Lecture Notes
  in Computational Science and Engineering}, vol.~45, pp. 353--356. Springer,
  Berlin, Germany (2005)

\bibitem{BennerSaak2013}
Benner, P., Saak, J.: Numerical solution of large and sparse continuous time
  algebraic matrix {R}iccati and {L}yapunov equations: a state of the art
  survey.
\newblock GAMM-Mitt. \textbf{36}(1), 32--52 (2013).
\newblock \doi{10.1002/gamm.201310003}

\bibitem{BlanesCasas2005}
Blanes, S., Casas, F.: On the necessity of negative coefficients for operator
  splitting schemes of order higher than two.
\newblock Appl.\ Numer.\ Math. \textbf{54}(1), 23--37 (2005).
\newblock \doi{10.1016/j.apnum.2004.10.005}

\bibitem{Caliari_etal2014}
Caliari, M., Kandolf, P., Ostermann, A., Rainer, S.: Comparison of software for
  computing the action of the matrix exponential.
\newblock BIT \textbf{54}(1), 113--128 (2014).
\newblock \doi{10.1007/s10543-013-0446-0}

\bibitem{DeLeo_etal2016}
De~Leo, M., Rial, D., de~la Vega, C.S.: High-order time-splitting methods for
  irreversible equations.
\newblock IMA J.~Numer.\ Anal. \textbf{36}(4), 1842--1866 (2016).
\newblock \doi{10.1093/imanum/drv058}

\bibitem{DruskinKnizhermanSimoncini2011}
Druskin, V., Knizhnerman, L., Simoncini, V.: Analysis of the rational {K}rylov
  subspace and {ADI} methods for solving the {L}yapunov equation.
\newblock SIAM J.~Numer.\ Anal. \textbf{49}(5), 1875--1898 (2011).
\newblock \doi{10.1137/100813257}

\bibitem{Guldogan_etal2016}
{G{\"u}ldo{\u g}an}, Y., {Hached}, M., {Jbilou}, K., {Kurulay}, M.: {Low rank
  approximate solutions to large-scale differential matrix Riccati equations}.
\newblock ArXiv e-prints  (2016).
\newblock \url{https://arxiv.org/abs/1612.00499}

\bibitem{Gustafsson1991}
Gustafsson, K.: Control-theoretic techniques for stepsize selection in explicit
  {R}unge-{K}utta methods.
\newblock ACM Trans.\ Math.\ Software \textbf{17}(4), 533--554 (1991).
\newblock \doi{10.1145/210232.210242}

\bibitem{Gustafsson_etal1988}
Gustafsson, K., Lundh, M., S{\"o}derlind, G.: A {PI} stepsize control for the
  numerical solution of ordinary differential equations.
\newblock BIT \textbf{28}(2), 270--287 (1988).
\newblock \doi{10.1007/BF01934091}

\bibitem{Hager1989}
Hager, W.W.: Updating the inverse of a matrix.
\newblock SIAM Rev. \textbf{31}(2), 221--239 (1989)

\bibitem{HairerWannerII}
Hairer, E., Wanner, G.: Solving ordinary differential equations. {II},
  \emph{Springer Series in Computational Mathematics}, vol.~14, second edn.
\newblock Springer, Berlin (1996).
\newblock \doi{10.1007/978-3-642-05221-7}

\bibitem{HansenOstermann2009_2}
Hansen, E., Ostermann, A.: High order splitting methods for analytic semigroups
  exist.
\newblock BIT \textbf{49}(3), 527--542 (2009).
\newblock \doi{10.1007/s10543-009-0236-x}

\bibitem{HeyouniJbilou2009}
Heyouni, M., Jbilou, K.: An extended block {A}rnoldi algorithm for large-scale
  solutions of the continuous-time algebraic {R}iccati equation.
\newblock Electron.\ Trans.\ Numer.\ Anal. \textbf{33}, 53--62 (2008/09)

\bibitem{HundsdorferVerwer2003}
Hundsdorfer, W., Verwer, J.: Numerical solution of time-dependent
  advection-diffusion-reaction equations, \emph{Springer Series in
  Computational Mathematics}, vol.~33.
\newblock Springer, Berlin (2003).
\newblock \doi{10.1007/978-3-662-09017-6}

\bibitem{IchikawaKatayama1999}
Ichikawa, A., Katayama, H.: Remarks on the time-varying {$H\sb \infty$}
  {R}iccati equations.
\newblock Syst.\ Control Lett. \textbf{37}(5), 335--345 (1999)

\bibitem{KoskelaMena2017}
{Koskela}, A., {Mena}, H.: {A Structure Preserving Krylov Subspace Method for
  Large Scale Differential Riccati Equations}.
\newblock ArXiv e-prints  (2017).
\newblock \url{https://arxiv.org/abs/1705.07507}

\bibitem{LangMenaSaak2015}
Lang, N., Mena, H., Saak, J.: On the benefits of the {$LDL^T$} factorization
  for large-scale differential matrix equation solvers.
\newblock Linear Algebra Appl. \textbf{480}, 44--71 (2015).
\newblock \doi{10.1016/j.laa.2015.04.006}

\bibitem{YidingSimoncini2015}
Lin, Y., Simoncini, V.: A new subspace iteration method for the algebraic
  {R}iccati equation.
\newblock Numer.\ Linear Algebra Appl. \textbf{22}(1), 26--47 (2015).
\newblock \doi{10.1002/nla.1936}

\bibitem{Mena_etal2017}
{Mena}, H., {Ostermann}, A., {Pfurtscheller}, L., {Piazzola}, C.: {Numerical
  low-rank approximation of matrix differential equations}.
\newblock ArXiv e-prints  (2017).
\newblock \url{https://arxiv.org/abs/1705.10175}

\bibitem{Petersen_etal2000}
Petersen, I.R., Ugrinovskii, V.A., Savkin, A.V.: Robust Control Design Using
  $H^{\infty}$ Methods.
\newblock Springer, London, UK (2000)

\bibitem{Saak2003}
Saak, J.: Effiziente numerische l\"{o}sung eines optimalsteuerungsproblems
  f\"{u}r die abk\"{u}hlung von stahlprofilen.
\newblock Master's thesis, Univ. Bremen, Bremen, Germany (2003)

\bibitem{MMESS101}
Saak, J., K\"{o}hler, M., Benner, P.: {M-M.E.S.S.-1.0.1} -- {T}he {M}atrix
  {E}quations {S}parse {S}olvers library.
\newblock DOI:10.5281/zenodo.50575 (2016).
\newblock See also: \url{www.mpi-magdeburg.mpg.de/projects/mess}

\bibitem{Simoncini2016}
Simoncini, V.: Computational methods for linear matrix equations.
\newblock SIAM Rev. \textbf{58}(3), 377--441 (2016).
\newblock \doi{10.1137/130912839}

\bibitem{SimonciniSzyldMonsalve2014}
Simoncini, V., Szyld, D.B., Monsalve, M.: On two numerical methods for the
  solution of large-scale algebraic {R}iccati equations.
\newblock IMA J.~Numer.\ Anal. \textbf{34}(3), 904--920 (2014).
\newblock \doi{10.1093/imanum/drt015}

\bibitem{Soderlind2002}
S{\"o}derlind, G.: Automatic control and adaptive time-stepping.
\newblock Numer.\ Algorithms \textbf{31}(1-4), 281--310 (2002).
\newblock \doi{10.1023/A:1021160023092}.
\newblock Numerical methods for ordinary differential equations (Auckland,
  2001)

\bibitem{SorensenZhou2002}
Sorensen, D.C., Zhou, Y.: Bounds on eigenvalue decay rates and sensitivity of
  solutions to {L}yapunov equations.
\newblock Tech. Rep. 02-07, Dept.\ of Comp.\ Appl.\ Math., Rice Univ., Houston,
  TX (2002).
\newblock \urlprefix\url{http://www.caam.rice.edu/caam/trs/tr02.html\#TR02-07}

\bibitem{Stillfjord2015}
Stillfjord, T.: Low-rank second-order splitting of large-scale differential
  {R}iccati equations.
\newblock IEEE Trans.\ Automat.\ Control \textbf{60}(10), 2791--2796 (2015).
\newblock \doi{10.1109/TAC.2015.2398889}

\end{thebibliography}
\end{document}